\documentclass[12pt,reqno]{amsart}
\usepackage{amsmath}
\usepackage{amsxtra}
\usepackage{amscd}
\usepackage{amsthm}
\usepackage{amsfonts}
\usepackage{amssymb}
\usepackage{eucal}

\theoremstyle{definition}

\theoremstyle{remark}

\theoremstyle{remark}

\numberwithin{equation}{section}

\let\Maketitle\maketitle
\def\maketitle{\hrule height0pt
\Maketitle\thispagestyle{empty}\let\maketitle\empty}

\def\fratop{\genfrac{}{}{0pt}1}

\let\Hat\widehat

\let\le\leqslant
\let\leq\leqslant
\let\ge\geqslant
\let\geq\geqslant

\def\>{\relax\ifmmode\mskip.666667\thinmuskip\relax\else\kern.111111em\fi}
\def\<{\relax\ifmmode\mskip-.333333\thinmuskip\relax\else\kern-.0555556em\fi}
\def\id{\mathrm{id}}
\def\vsk#1>{\vskip#1\baselineskip}
\def\vv#1>{\vadjust{\vsk#1>}\ignorespaces}
\def\vvn#1>{\vadjust{\nobreak\vsk#1>\nobreak}\ignorespaces}
\def\vvgood{\vadjust{\penalty-500}}

\def\qdet{\mathop{\mathrm{qdet}}\nolimits}

\def\XXX/{{\slshape XXX\/}}
\def\xxx{\hbox{\tiny\sl XXX\/}}
\def\Hxxx{H_{\xxx}}
\let\alb\allowbreak
\let\vp\vphantom
\def\tbigcap{\mathrel{\textstyle\bigcap}}

\makeatletter
\def\section{\def\@secnumfont{\mdseries}\@startsection{section}{1}%
 \z@{.7\linespacing\@plus\linespacing}{.5\linespacing}%
 {\normalfont\scshape\centering}}
\def\subsection{\def\@secnumfont{\bfseries}\@startsection{subsection}{2}%
 \z@{.5\linespacing\@plus.7\linespacing}{-.5em}%
 {\normalfont\bfseries}}
\makeatother

\let\Smallskip\smallskip
\let\Medskip\medskip
\let\Bigskip\bigskip
\def\smallskip{\par\Smallskip}
\def\medskip{\par\Medskip}
\def\bigskip{\par\Bigskip}

\let\nc\newcommand

\nc{\wt}{\widetilde}
\nc{\on}{\operatorname}
\nc{\ch}{\mbox{ch}}
\nc{\Z}{\mathbb Z}
\nc{\C}{\mathbb C}
\nc{\R}{\mathbb R}
\nc{\pone}{\mathbb{CP}^1}
\nc{\p}{\partial}
\nc{\Bc}{\mathcal B}
\nc{\F}{\mathcal F}
\nc{\arr}{\rightarrow}
\nc{\larr}{\longrightarrow}
\nc{\al}{\alpha}
\nc{\ri}{\rangle}
\nc{\lef}{\langle}
\nc{\la}{\lambda}
\nc{\ep}{\epsilon}
\nc{\eps}{\varepsilon}
\nc{\Om}{\Omega}
\nc{\om}{\omega}
\nc{\La}{\Lambda}
\nc{\Gm}{\Gamma}
\nc{\gm}{\gamma}
\nc{\si}{\sigma}
\nc{\sh}{\sigma}

\nc{\D}{\mathcal D}
\nc{\A}{{\mathbb A}}
\nc{\g}{\mathfrak g}
\nc{\h}{\mathfrak h}
\nc{\mg}{\mathfrak m}
\nc{\n}{\mathfrak n}
\nc{\z}{\mathfrak Z}
\nc{\mc}{\mathcal}
\nc{\M}{\mathcal M}
\nc{\N}{\widehat\n}
\nc{\G}{\widehat\g}
\nc{\De}{\Delta_+}
\nc{\gt}{\widetilde{\g}}
\nc{\one}{\mathbf 1}
\nc{\Hh}{\mathcal H_\beta}
\nc{\qp}{q^{\frac{k}{2}}}
\nc{\qm}{q^{-\frac{k}{2}}}
\nc{\qn}{\frac{[m]_q^2}{[2m]_q}}
\nc{\cri}{_{\on{cr}}}
\nc{\kk}{h^\vee}
\nc{\sun}{\widehat{\sw}_N}
\nc{\hh}{\widehat{\mathfrak h}}
\nc{\HH}{{\mathcal H}_{q,t}}
\nc{\ca}{\wt{{\mathcal A}}_{h,k}(\sw_2)}
\nc{\Res}{\on{Res}}

\nc{\bean}{\begin{eqnarray}}
\nc{\eean}{\end{eqnarray}}
\nc{\be}{\begin{equation*}}
\nc{\ee}{\end{equation*}}
\nc{\bea}{\begin{eqnarray*}}
\nc{\eea}{\end{eqnarray*}}
\nc{\beq}{\begin{equation}}
\nc{\eeq}{\end{equation}}
\nc{\bs}{\boldsymbol}
\nc{\Ref}[1]{{\rm(\ref{#1})}}

\nc{\glN}{\mathfrak{gl}_N}
\nc{\glt}{\mathfrak{gl}_2}

\nc{\End}{{\rm End}}
\nc{\T}{{\bs{p}}}
\nc{\nash}{M_{\bs \La}[\La^{(\infty)}]}
\nc{\Sing}{{\rm Sing}}
\nc{\snash}{\Sing\,\W[\>l\>]}
\nc{\Snash}{\Sing\,L_{\bs \La}[\>l\>]}
\nc{\Y}{Y(\glt)}
\nc{\Ys}{Y(\mathfrak{sl}_2)}
\nc{\Cx}{\C[x_1,\dots,x_n]}
\nc{\MLz}{M_{\bs\La}(\bs z)}
\nc{\Xe}{\mathcal X}

\nc{\Tee}{\mathcal T}
\nc{\Te}{\mathfrak T}
\nc{\Se}{\mathfrak S}

\nc{\AD}{A_{\T,D}}
\nc{\AT}{A_{\T,T}}
\nc{\Dh}{\mathcal{D}_{\bs h}}
\nc{\Dhp}{\mathcal{D}_{\bs h(\bs p)}}

\nc{\qa}{p(u,\bs a)}
\nc{\qap}{p(u,\bs a(\bs p))}

\def\Po{P}
\nc{\DD}{\mathfrak{D}}
\nc{\Bg}{\mathfrak{B}}

\nc{\tT}{\tilde{T}}
\nc{\Cl}{{\C_l[y_1,\dots,y_{n-1}]^{{\rm Sym}}}}
\nc{\W}{W_{a,d}}
\nc{\Wy}{\mathcal {W}_{a,d}}
\nc{\qat}{\tilde {p}(u,\tilde {\bs a})}
\nc{\Wr}{{\rm Wr}}

\nc{\SC}{\Sing\>(\C^2)^{\otimes n}[\>l\>]}
\nc{\SCf}{\Sing\>(\C^2)^{\otimes 4}[\>2\>]}
\nc{\bv}{\om(t_1,\dots,t_l)}
\nc{\SSnash}{\Sing\,L_{\bs 1}[\>l\>]}

\let\shift\vartheta

\begin{document}

\title[Bethe Algebra of Homogeneous \XXX/ Model Has Simple Spectrum]
{Bethe Algebra of Homogeneous \XXX/ Heisenberg Model Has Simple
Spectrum}

\author[E.\,Mukhin, V.\,Tarasov, and A.\>Varchenko]
{E.\,Mukhin ${}^{*,1}$, V.\,Tarasov ${}^{*,\star,2}$,
\and A.\>Varchenko {${}^{**,3}$} }

\thanks{${}^1$\ Supported in part by NSF grant DMS-0601005}
\thanks{${}^2$\ Supported in part by RFFI grant 08-01-00638}
\thanks{${}^3$\ Supported in part by NSF grant DMS-0555327}

\maketitle

\centerline{\it ${}^*$Department of Mathematical Sciences,
Indiana University -- Purdue University,}
\centerline{\it Indianapolis, 402 North Blackford St, Indianapolis,
IN 46202-3216, USA\/}
\smallskip
\centerline{\it $^\star$St.\,Petersburg Branch of Steklov Mathematical
Institute}
\centerline{\it Fontanka 27, St.\,Petersburg, 191023, Russia\/}
\smallskip
\centerline{\it ${}^{**}$Department of Mathematics, University of
North Carolina at Chapel Hill,} \centerline{\it Chapel Hill, NC
27599-3250, USA\/}

\bigskip
\begin{abstract}
We show that the algebra of commuting Hamiltonians of the
homogeneous \XXX/ Heisenberg model has simple spectrum on the subspace
of singular vectors of the tensor product of two-dimensional $\glt$-modules.
As a byproduct we show that there exist exactly
$\binom {n}{l}-\binom{n}{l-1}$ two-dimensional vector subspaces
$V \subset \C[u]$ with a basis $f,g\in V$ such that
$\deg f = l, \deg g = n-l+1$ and $f(u)g(u-1) - f(u-1)g(u) = (u+1)^n$.
\end{abstract}

\section{Introduction}
\subsection{Homogeneous \XXX/ Heisenberg model}
Consider the vector space $(\C^2)^{\otimes n}$ and the linear operator
\vvn-.5>
\be
\Hxxx\,=\,-\sum_{j=1}^n\,(\>\sh_1^{(j)}\sh_1^{(j+1)}+\>
\sh_2^{(j)}\sh_2^{(j+1)}+\>\sh_3^{(j)}\sh_3^{(j+1)}\>)\;,
\ee
where \;$\sh_a^{(k)}=1^{\otimes(k-1)}\otimes\sh_a\otimes1^{(n-k)}$,
\;$\sh_a^{(n+1)}=\sh_a^{(1)}$, \;and \;$\sh_1,\,\sh_2,\,\sh_3$ \,are
\vvn.2>
the Pauli matrices,
\be
\sh_1\,=\begin{pmatrix}\,0&1\,\\\,1&0\,\end{pmatrix}\;,\qquad
\sh_2\,=\begin{pmatrix}\,0&-i\>\\\,i&0\>\end{pmatrix}\;,\qquad
\sh_3\,=\begin{pmatrix}\,1&0\\\,0&-1\end{pmatrix}\;.
\vv.2>
\ee
The operator $\Hxxx$ is the Hamiltonian of the celebrated \XXX/ Heisenberg
model, also called the homogeneous \XXX/ model, and the problem is to find
eigenvalues and eigenvectors of the Hamiltonian.

\smallskip
This problem was first addressed in the pioneering work \cite{Bt} by H.\,Bethe,
who looked for eigenvectors of $\Hxxx$ in a certain special form. His method
\vvgood
and its further extensions are traditionally called the Bethe ansatz.
The current literature on the \XXX/ model and its generalizations, {\sl XXZ\/}
and {\sl XYZ\/} models, as well as their counterparts in statistical mechanics,
the six- and eight-vertex models, is enormous. We limit ourselves to mentioning
just two books, \cite{Ba} and \cite{KBI}. However, even numerous references
therein hardly cover a half of the bibliography on the subject.

\smallskip
The Hamiltonian $\Hxxx$ can be included into a one\>-parameter family of
commuting linear operators called the transfer matrix, see~\cite{Ba},
\cite{FT}, \cite{KBI}. We call a commutative unital subalgebra of linear
operators on $(\C^2)^{\otimes n}$ generated by the transfer matrix the Bethe
algebra. The actual problem is to construct eigenvalues and eigenvectors for
the Bethe algebra.

The elements of the Bethe algebra commute with the natural $\glt$-action on
$(\C^2)^{\otimes n}$. Therefore, the eigenspaces of the Bethe algebra are
representations of $\glt$, and it suffices to construct highest weight vectors
of those representations.

\smallskip
The Bethe ansatz method associates to every admissible solution
$(\la_1,\dots\la_l)$ of the system of equations
\vvn-.3>
\beq
\label{bae}
\left(\frac{\la_j+\frac i2}{\la_j-\frac i2}\>\right)^{\!n}\,=\,
\prod_{\fratop{k=1}{k\ne j}}^l\,\frac{\la_j-\la_k+i}{\la_j-\la_k-i}\;,
\qquad j=1,\dots,l\;,
\eeq
a vector in $(\C^2)^{\otimes n}$, called the corresponding Bethe vector,
see~\cite{FT}. A solution $(\la_1,\dots,\la_l)$ is called admissible if all
$\la_1,\dots,\la_l$ are distinct, and all factors in \Ref{bae} are nonzero.
A nonzero Bethe vector
is a highest weight vector of an $(n-2l+1)$-dimensional irreducible
representation of $\glt$, and all vectors in that representation are
eigenvectors of each element of the Bethe algebra sharing the same eigenvalue.

\smallskip
It is an important question whether the Bethe ansatz method produces all
eigenvectors of the Bethe algebra. This question is referred to as the question
of completeness of the Bethe ansatz for finite chains. It was discussed by
H.\,Bethe himself in~\cite{Bt} and many times since then by other authors.
For instance, see a recent discussion in~\cite{B2}. However, no rigorous proof
is available even for the so-called inhomogeneous models. Moreover, as one
can see from the results of this paper, Sklyanin's separation of variables does
not prove completeness of the Bethe ansatz to the very end, though it is indeed
an important step towards the proof.

To be more precise, there are certain quantum integrable models for which
the completeness of the Bethe ansatz has been proved. For example,
see~\cite{YY} and Theorem~1.2.2 in~\cite{KBI}. The proofs for those models
are based on a variational principle and convexity of some auxiliary action.
However, for the the \XXX/ model, the corresponding action is not convex,
and that technique fails.

\smallskip
In this paper we establish the completeness of the Bethe ansatz method for
the homogeneous \XXX/ model provided the method is improved in a certain way,
see below in the introduction. We show that the spectrum of the Bethe algebra
of the homogeneous \XXX/ model is simple, that is, all eigenspaces of the Bethe
algebra are irreducible $\glt$-modules. We also show that eigenvalues of
the Bethe algebra are in a one\>-to\>-one correspondence with certain
second-order linear difference equations with two linearly independent
polynomial solutions. We prove similar results for inhomogeneous higher spin
\XXX/ models.

\smallskip
To continue with an introduction and match the notation in the main part
of the paper, we change the variables in system~\Ref{bae},
\be
\la_j\,=\,\frac i2\,(t_j+1)\,,\qquad j=1,\dots,l\,,
\ee
and write the system in the polynomial form
\beq
\label{BAEi}
(t_j+2)^n\prod_{\fratop{k=1}{k\ne j}}^l\,(t_j - t_k - 1)\,=\,
(t_j+1)^n\prod_{\fratop{k=1}{k\ne j}}^l\,(t_j - t_k + 1)\;,
\qquad j = 1,\dots,l\;.
\eeq
We call system~\Ref{BAEi} the system of the Bethe ansatz equations. The system
is invariant with respect to permutations of $t_1,\dots,t_l$, so the symmetric
group $S_l$ acts on solutions to the Bethe ansatz equations.

We denote by $\om(t_1,\dots,t_l)$ the Bethe vector corresponding to
an admissible solution $(t_1,\dots, t_l)$ of the Bethe ansatz equations.
The Bethe vectors corresponding to admissible solutions with permuted
coordinates are equal. The number of Bethe vectors $\om(t_1,\dots,t_l)$ is
equal to the number of $S_l$-orbits of admissible solutions to
system~\Ref{BAEi}.

Since each element of the Bethe algebra commutes with the natural $\glt$-action
on $(\C^2)^{\otimes n}$, it is enough to diagonalize the action of the Bethe
algebra on each subspace of $\glt$-singular vectors of given weight,
\be
\SC\,=\,\{\,v\in (\C^2)^{\otimes n}\ | \ \,e_{12}v=0, \ \,e_{11}v=(n-l)\>v,
\ \,e_{22} v = l\>v\,\}\;,
\ee
with \,$2l\le n$\>. For every admissible solution $(t_1,\dots, t_l)$ of
the Bethe ansatz equations, the Bethe vector $\om(t_1,\dots,t_l)$ belongs to
the subspace \,$\SC$\,.

\smallskip
To illustrate the problem with completeness of the Bethe ansatz in the standard
form and the way it can be resolved, let us consider an example.

\smallskip
Let $n=4$ and $l=2$. Then \,$\dim\,\SCf=2$\>, the operator $\Hxxx$ restricted
to \,$\SCf$ \,has eigenvalues \>$5$ \>and \>$-3$\>. The Bethe ansatz equations
are
\begin{align}
\label{BAEii}
(t_1+2)^4\,(t_1 - t_2 - 1)\, &{}=\,(t_1+1)^4\,(t_1 - t_2 + 1)\;,
\\[3pt]
(t_2+2)^4\,(t_2 - t_1 - 1)\, &{}=\,(t_2+1)^4\,(t_2 - t_1 + 1)\;,
\notag
\end{align}
and there is only one orbit of admissible solutions:
\beq
\label{tt12}
t_1\,=\,-\>\frac 32\>+\>\frac 12\>\sqrt{-\>\frac 13}\ ,\qquad
t_2\,=\,-\>\frac 32\>-\>\frac 12\>\sqrt{-\>\frac 13}\ .
\eeq
The Bethe vector $\om(t_1,t_2)$ is an eigenvector of $\Hxxx$ with eigenvalue
$5$.

\smallskip
The results of this paper say that each eigenspace of $\Hxxx$ acting on $\SCf$
corresponds to a difference equation
\beq
\label{deq4}
u^4 f(u)\>-\>\Bc(u)f(u-1)\>+\>(u+1)^4f(u-2)\,=\,0\,,
\eeq
where $\Bc(u)$ is a polynomial, and the difference equation has polynomial
solutions of degree $2$ and $3$. The corresponding eigenvalue of $\Hxxx$ equals
\,$1-2\>\Bc'(0)/\Bc(0)$\>.
\vsk.1>
Indeed, there are exactly two such difference equations. The first one has
$\Bc(u)\>=\>2u^4+4u^3-2u+1$, and solutions \,$u^2+3u+\frac73$ \,and
\,$u^3+6u^2+11u+\frac{13}2$\,. The roots of the quadratic polynomial
are numbers $t_1$ and $t_2$ given by~\Ref{tt12}.

The second difference equation~\Ref{deq4} with polynomial solutions of degree
$2$ and $3$ has $\Bc(u)\>=\>2u^4+4u^3-2u-1$, and solutions \,$(u+1)\>(u+2)$
\,and \,$u^3+6u^2+10u+\frac92$\,. The roots of the quadratic polynomial,
$t_1=-1$ and $t_2=-2$, form a nonadmissible solution to system~\Ref{BAEii}\>,
and the Bethe vector $\om(t_1,t_2)$ for $t_1=-1$, $t_2=-2$ equals zero.

\smallskip
For general $n$ and $l$ such that $2l\le n$, the results of this paper for
the homogeneous \XXX/ model say that eigenspaces of the Bethe algebra acting
on $\SC$ are one-dimensional. They are in a one\>-to\>-one correspondence
with difference equations
\vvn.3>
\beq
\label{deq}
u^n f(u)\>-\>\Bc(u)f(u-1)\>+\>(u+1)^nf(u-2)\,=\,0\,,
\vv.3>
\eeq
where $\Bc(u)$ is a polynomial, and those difference equations have polynomial
solutions of degree $l$ and $n-l+1$. The corresponding eigenvalues of elements
of the Bethe algebra are described by the polynomial $\Bc(u)$. In particular,
the eigenvalue of $\Hxxx$ equals \,$1-2\>\Bc'(0)/\Bc(0)$\>. The roots
$t_1,\dots t_l$ of the polynomial solution of equation~\Ref{deq} of degree $l$
form a solution of system~\Ref{BAEi}. The Bethe vector \,$\om(t_1,\dots t_l)$
\,is nonzero if and only if the solution $(t_1,\dots t_l)$ is admissible.

\smallskip
To obtain an eigenvector of the Bethe algebra corresponding to a difference
equation~\Ref{deq} with two polynomial solutions, we use the following
construction. The space $(\C^2)^{\otimes n}$ has a structure of a module over
the Yangian $\Y$, and the Bethe algebra of the homogeneous \XXX/ model is
the image of a commutative subalgebra \,$\Bg\subset\Y$, called the Bethe
subalgebra.
We take another $\Y$-module $\W$, described in Section~\ref{holomorphic},
which is the holomorphic representation of $\Y$ associated with the polynomials
$a(u)=(u+1)^n$ and $d(u)=u^n$. There is a natural epimorphism
$\sh:\W\to(\C^2)^{\otimes n}$ of $\Y$-modules.

Using the roots $t_1,\dots,t_l$ of the polynomial solution of
equation~\Ref{deq} of degree $l$ and Sklyanin's procedure of separation of
variables~\cite{Sk}, we define a nonzero vector \,$\tilde\omega(t_1,\dots,t_l)$
in $\W$, which is an eigenvector of \,$\Bg$ \>acting on $\W$. We consider
the maximal \,$\Bg$-invariant subspace $V\subset\W$ that contains
\,$\tilde\omega(t_1,\dots,t_l)$ \>and does not contain other linearly
independent eigenvectors of \,$\Bg$. We show that the image
$\sh(V)\subset(\C^2)^{\otimes n}$ is a one-dimensional subspace of $\SC$.
Since \,$\sh$ \>is an homomorphism of $\Y$-modules, $\sh(V)$ is an eigenspace
of the Bethe algebra acting on $\SC$ with the same eigenvalues as the
eigenvalues of $\tilde\omega(t_1,\dots,t_l)$ with respect to the action
of $\Bg$ on $\W$. The subspace $\sh(V)\subset\SC$ is that one-dimensional
subspace of eigenvectors which we assigned to difference equation~\Ref{deq}
with two polynomial solutions.

If $(t_1,\dots,t_l)$ is an admissible solution, then the subspace $V\subset\W$
is one-dimensional, and the subspace $\sh(V)$ is spanned by the Bethe vector
\,$\omega(t_1,\dots,t_l)$.

\smallskip
The construction described above provides a generalization of the Bethe ansatz
method in which the solutions to the Bethe ansatz equations are replaced by
difference equation~\Ref{deq} with two polynomial solutions, and the Bethe
vectors in $\SC$ are replaced by the subspaces $\sh(V)$. Our result says that
the generalized Bethe vectors form a basis in $\SC$ and, moreover, the spectrum
of the Bethe algebra is simple.

\smallskip
As a remark, we would like to indicate another way to obtain the eigenspace
of the Bethe algebra acting on $\SC$ corresponding to the difference
equation~\Ref{deq}. We may consider the inhomogeneous \XXX/ model depending
on parameters $z_1,\dots,z_n$. The corresponding system of the Bethe ansatz
equations are
\beq
\label{BAEz}
\prod_{s=1}^n\,(t_j-z_s+2)\,
\prod_{\fratop{k=1}{k\ne j}}^l\,(t_j - t_k - 1)\,=\,
\prod_{s=1}^n\,(t_j-z_s+1)\,
\prod_{\fratop{k=1}{k\ne j}}^l\,(t_j - t_k + 1)\;,
\eeq
$j=1,\dots,l$\>. It follows from the results of this paper that if
$f(u)=\prod_{j=1}^l(u-t_j)$ is a solution of the difference equation~\Ref{deq},
then for generic $\bs z=(z_1,\dots,z_n)$ there exists an admissible solution
\,$\bs t(\bs z)=(t_1(\bs z),\dots,t_l(\bs z))$ of system~\Ref{BAEz} such that
\,$\bs t(\bs z)\to(t_1,\dots t_l)$ as $\bs z\to 0$. The Bethe vectors
\,$\om(\bs t(\bs z);\bs z)$ are nonzero for generic $\bs z$, and
the eigenspaces \,$\C\>\om(\bs t(\bs z);\bs z)$ have a one-dimensional limit as
$\bs z\to 0$, which is the eigenspace of the Bethe algebra of the homogeneous
\XXX/ model. A similar approach for the $\glN$ Gaudin model is developed
in~\cite{MTV5}.

\smallskip
The correspondence between the eigenvectors of the Bethe algebra and
second-order linear difference equation with two polynomial solutions is in
the spirit of the geometric Langlands correspondence in which eigenfunctions
of commuting differential operators correspond to connections on curves.

\medskip
Equation~\Ref{deq} is known in the physical literature as Baxter's equation.
Its connection with the Bethe ansatz equations has been studied in many papers.
The fact that the roots of a polynomial solution of Baxter's equation give
a solution of the Bethe ansatz equations (provided the roots are distinct)
is known as Manakov's principle and the analytic Bethe ansatz. An important
observation about the existence of a second polynomial solution of Baxter's
equation has been done in~\cite{PS}. A similar observation in a much more
general context has been made independently in~\cite{MV2}, \cite{MV3}.

\subsection{Content of the paper}
The results of this paper for the \XXX/ model are discrete analogues of
the results of \cite{MTV3} for the Gaudin model.

\smallskip
In Section \ref{sec yangian} we discuss the Yangian $\Y$, the Bethe subalgebra
$\Bg\subset\Y$, and Yangian modules. In particular, we describe the holomorphic
representation $\W$ of the Yangian $\Y$. The module $\W$ is associated with
two monic polynomials
\be
a(u)\,=\,\prod_{i=1}^n\,(u-z_i+m_i)\qquad {\rm and}
\qquad d(u)\,=\,\prod_{i=1}^n\,(u-z_i)
\ee
and is isomorphic to $\C[x_1,\dots,x_n]$ as a vector space.

We introduce a collection $\bigl((m_1,0),\dots,(m_n,0)\bigr)$ of $\glt$-weights
and say that the pair $\bigl((m_1,0), \dots, (m_n,0)\bigr),\,l$ \,is separating
if $\sum_{i=1}^n m_i - 2l +1 + s\,\ne\,0$ for all $s=1,\dots,l$.

\medskip
In Sections~\ref{Universal difference operator}\,--\>\ref{Multiplication Mult},
we study the algebras $A_W$ and $A_D$, and relations between them.
Eventually, we show that the algebras $A_W$ and $A_D$ are isomorphic,
see~Theorem~\ref{1st main thm}.

\smallskip
The algebra $A_W$ is the image of the Bethe subalgebra $\Bg$ acting on
the subspace $\snash\subset\W$ of $\glt$-singular vectors. We consider
a polynomial $B(u,\bs H)\>=\>2u^n + H_1u^{n-1}+\dots+H_n$, whose coefficients
$H_k\in\End\,(\Sing\,\W[\>l\>])$ are generators of $A_W$, and introduce
the universal difference operator
\vvn.2>
\be
\DD_{\>\snash}\,=\,d(u)\>-\>{B(u, \bs H)}\,\shift^{-1}+\>{a(u)}\,\shift^{-2}
\vv.3>
\ee
acting on $\snash$-valued functions in $u$. Here $\shift:f(u)\mapsto f(u+1)$.

\smallskip
The algebra $A_D$ is defined in Section \ref{Algebra A_D}. We consider
the space $\C^{l+n}$ with coordinates $\bs a=(a_1,\dots,a_l)$ and
$\bs h=(h_1,\dots,h_n)$, polynomials $B(u,\bs h)\,=\,2u^n+h_1u^{n-1}+\dots+h_n$
and $\qa=u^l+a_1u^{l-1}+\dots+a_l$, and the difference operator
\vvn.2>
\be
\D_{\bs h}\,=\,d(u)\>-\>{B(u,\bs h)}\,\shift^{-1}+\>{a(u)}\,\shift^{-2}
\vv.3>
\ee
We define the scheme $C_D $ of points $\T\in \C^{l+n}$ such that the polynomial
$p(u,\bs a(\T))$ lies in the kernel of the difference operator
$\D_{\bs h(\T)}$. The algebra $A_D$ is the algebra of functions on $C_D$.
There is a natural epimorphism \,$\psi_{DW}:A_D \to A_W$ \,such that
\,$\psi_{DW}(h_k)=H_k$, \,see Theorem~\ref{thm D to B}.

\smallskip
Using the Bethe ansatz method, we prove that if $z_1,\dots,z_n$ are generic
and the pair $\bigl((m_1,0),\dots,(m_n,0)\bigr),\,l$ \,is separating, then the
scheme $C_D$ considered as a set has at least \,$\dim\,\snash$ distinct points,
see Section~\ref{Bethe ansatz and C_D}.

\smallskip
In Section \ref{On separation for slt}, we review Sklyanin's procedure of
separation of variables in the \XXX/ model and construct the universal weight
function. Theorem \ref{thm on Bethe ansatz} connects the algebras $A_D,\,A_W$
and the universal weight function.

\smallskip
The algebra $A_D$ acts on itself by multiplication operators. We denote by
$L_f$ the operator of multiplication by an element $f\in A_D$. The algebra
$A_D$ acts on its dual space $A_D^*$ by operators $L^*_f$\>,
dual to multiplication operators.
Using the universal weight function
we define a linear map \,$\tau:A_D^*\to\snash$ \,and prove that if the pair
$\bigl((m_1,0),\dots,(m_n,0)\bigr),\,l$ \,is separating, then $\tau$ is
an isomorphism that intertwines the action of operators $L_f^*$\>,
\;$f\in A_D$, with the action of operators $\psi_{DW}(f)\in\End(\snash)$,
see Theorem~\ref{1st main thm}. Therefore, we prove that
$\psi_{DW}:A_D\to A_W$ is an algebra isomorphism.
Theorem~\ref{1st main thm} is our first main result.

Using the Grothendieck residue, we define an isomorphism $\phi:A_D\to A_D^*$
of $A_D$-modules, see Section~\ref{Grothendieck bilinear form on A-D}.
Therefore, if the pair $\bigl((m_1,0),\dots,(m_n,0)\bigr),\,l$ \,is separating,
the composition $\tau\phi\,:\,A_D\to\snash$ is a linear isomorphism which
intertwines the action of the algebra $A_D$ on itself by multiplication
operators and the action of the Bethe algebra $A_W$ on $\snash$.

\medskip
In Sections \ref{Three more algebras}
through~\ref{Equations with polynomial solutions only},
we impose more conditions on $m_1,\dots,m_n$ and
$z_1,\dots,z_n$. We assume that $m_1,\dots,m_n$ are natural numbers.
We keep the assumption that the pair $\bigl((m_1,0),\dots,(m_n,0)\bigr),\,l$
is separating, that takes the form $2l\le\sum_{s=1}^n m_s$\>.
We also assume that $z_i-z_j\notin \Z$ if $i\ne j$.

\smallskip
In Sections~\ref{Three more algebras}\,--\,\ref{Equations with polynomial
solutions only}, we study three more algebras $A_G$, \,$A_P$ and $A_L$,
and relations between them. The algebra $A_G$ is defined in
Section~\ref{Three more algebras}. We consider the subspace
$\C_d[u]\subset\C[u]$ of all polynomials of degree $\leq d$ \,for a suitably
large number $d$, \>and the Grassmannian of all two-dimensional subspaces of
$\C_d[u]$. Using the numbers $z_1,\dots,z_n$ and $m_1,\dots,m_n$ we define
\vvn.1>
\>$n+1$ \,Schubert cycles $C_{\F(z_1),\La^{(1)}}\>,\dots,
\>C_{\F(z_n),\La^{(n)}}\>,\,C_{\F(\infty), \La^{(\infty)}}$
in the Grassmannian. The algebra $A_G$ is the algebra of functions
on the intersection of the Schubert cycles.

\smallskip
The algebra $A_P$ is defined in Section~\ref{sectAP}.
\vvn.2>
Let \,$\tilde l=\sum_{s=1}^n m_s+1-l$,
\;$\tilde{\bs a}=(\tilde a_1,\dots,\alb\tilde a_{\tilde l-l-1},
\tilde a_{\tilde l-l+1},\dots,\tilde a_{\tilde l})$\>,
\;$\qat=u^{\tilde l}+\tilde a_1u^{\tilde l-1}+\dots+
\tilde a_{\tilde l-l-1}u^{l+1}+\tilde a_{\tilde l-l+1}u^{l-1}+\dots+
\tilde a_{\tilde l}$\>,
\vvn.1>
\,and consider the space $\C^{\tilde l + l + n - 1}$
with coordinates $\tilde{\bs a}, \bs a, \bs h$. We define
the scheme $C_P$ as the scheme of points $\T \in \C^{\tilde l + l + n - 1}$
such that the polynomials $p(u,\bs a(\T))$ and $\tilde p(u,\tilde{\bs a}(\T))$
lie in the kernel of the difference operator $\D_{\bs h(\T)}$.
\vvn.1>
The algebra $A_P$ is the algebra of functions on $C_P$. The map
\vvn.1>
\,$(p(u,\bs a(\T)),\tilde p(u,\tilde{\bs a}(\T)),\D_{\bs h(\T)})\mapsto
(p(u,\bs a(\T)),\D_{\bs h(\T)})$ \,defines a natural epimorphism
$\psi_{DP}:A_D \to A_P$.
We also show that the algebras $A_G$ and $A_P$ are naturally isomorphic.

To define the algebra $A_L$, see Section~\ref{A_L},
we consider the tensor product
\vvn.2>
\be
L_{\bs \La}(\bs z)\,=\,L_{\La^{(1)}}(z_1)\otimes\dots\otimes L_{\La^{(n)}}(z_n)
\vv.2>
\ee
of evaluation Yangian modules, where \,$L_{\La^{(i)}}$ is the irreducible
$\glt$-module of highest weight $\La^{(i)}=(m_i,0)$\>. \,The algebra $A_L$
is the image of the Bethe subalgebra $\Bg\subset\Y$ acting on the subspace
$\Snash\subset L_{\bs \La}(\bs z)$ of $\glt$-singular vectors. The $\Y$-module
$L_{\bs \La}(\bs z)$ is isomorphic to the quotient module $\W/K$\>, where
$K\subset\W$ is the kernel of the Yangian Shapovalov form on $\W$.
We denote by $\sh : \snash \to \Snash$ the epimorphism of vector spaces
corresponding to the epimorphism $\W\to L_{\bs \La}(\bs z)$ of $\Y$-modules.
The epimorphism $\sh$ induces the algebra epimorphism $\psi_{WL}:A_W\to A_L$\>.

\smallskip
We denote by \,$\xi:A_D\to\Snash$ \,the composition of maps \,$\sh\tau\phi$,
and by \,$\psi_{DL}:A_D\to A_L$ \,the composition of maps $\psi_{WL}\psi_{DW}$.
We show that the kernels of the maps \>$\xi$\>, $\psi_{DL}$ \>and $\psi_{DP}$
coincide. This allows us to obtain an algebra isomorphism
$\psi_{PL}:A_P\to A_L$ and a linear isomorphism $\zeta:A_P\to\Snash$
intertwining the action of $A_P$ on itself by multiplication operators and the
action of the Bethe algebra $A_L$ on $\Snash$. This is our second main result,
see Theorem \ref{second main thm}.

In Section~\ref{Equations with polynomial solutions only}, we use the Yangian
Shapovalov form on $L_{\bs \La}(\bs z)$ and the map $\zeta$ to obtain
a linear isomorphism $\theta:A_P^*\to\Snash$ intertwining the action of
operators $L_f^*$\>, \;$f\in A_P$, with the action of operators
$\psi_{PL}(f)\in\End(\Snash)$, see Theorem~\ref{useful cor}.
Using the isomorphism, we show that eigenvectors of the action of
the algebra $A_L$ on $\Snash$ are in a one\>-to\>-one correspondence with
certain second-order linear difference equations with two polynomial solutions
of degrees $l$ and $n-l+1$, see Corollary~\ref{Cor 2 of 2nd main thm}.

\smallskip
Section \ref{Heisenberg chain} contains the analogues of the previous
results for the homogeneous \XXX/ Heisenberg model.

\smallskip
We recapitulate the main results of this paper as three commutative diagrams.
The horizontal arrows of the diagrams are isomorphisms, the downward vertical
arrows are epimorphisms, and the upward vertical arrow is an embedding.

The first diagram shows the algebras of functions $A_D$, $A_P$ on difference
operators with respectively one or two polynomials in the kernels, the algebra
$A_G$ of functions on the intersection of Schubert cycles, the Bethe algebras
$A_W$ and $A_L$, associated with $\snash$ and $\Snash$, respectively, and their
homomorphisms:
\be
\begin{CD}
@. A_D @>\psi_{DW}>> A_W\\
@. @V\psi_{DP}VV @VV\psi_{WL}V\\
A_G %\,\leftarrow\kern-1.5em
@>>\psi_{GP}> A_P @>>\psi_{PL}> A_L @.\kern4.8em
\end{CD}
\vv.3>
\ee
The other two diagram show the vector spaces involved:
\vvn.3>
\be
\begin{CD}
A_D^* @>\tau>>\,\snash\,\\
@A(\psi_{DP})^*AA @VV\sh V\\
A_P^* @>>\theta> \Snash
\end{CD}
\kern6.6em
\begin{CD}
A_D @>\tau\phi>>\,\snash\,\\
@V\psi_{DP} VV @VV\sh V\\
A_P @>>\zeta> \Snash
\end{CD}
\ee
Each vector space on these diagrams is a module over the corresponding algebra
on the first diagram, and all linear maps are consistent with the algebra
homomorphisms.

\subsection*{Acknowledgments}
The authors thank referees for helpful comments.

\section{Yangian $Y(\glt)$ and Yangian modules}
\label{sec yangian}

\subsection{Lie algebra $\glt$}
Let $e_{ab}$, $\;a,b=1,2$, be the standard generators of
the complex Lie algebra $\glt$.
We have \,$\glt=\n^+\oplus\h\oplus\n^-$ where
\vvn.2>
\be
\n^+\>=\,\ \C\cdot e_{12}\,, \qquad
\h\,=\,\C \cdot e_{11}\ \oplus \ \C \cdot e_{22}\, , \qquad
\n^-\>=\, \C \cdot e_{21}\,.
\ee

\smallskip
For a $\glt$-weight $\La\, \in \h^*$, we
denote by $M_\La$ the Verma $\glt$-module with highest weight $\La$
and by $L_\La$ the irreducible $\glt$-module with highest weight $\La$.

\subsubsection{}
Let $\bs \La = (\La^{(1)},\dots,\La^{(n)})$ be a collection of $\glt$-weights,
where $\La^{(i)}=(\La_1^{(i)},\La_2^{(i)})$ for \,$i=1,\dots, n$.
Let $l$ be a nonnegative integer. The pair $\bs \La$, $l$ will be called
{\it separating\/} if \ $\sum_{i=1}^n(\La_1^{(i)}-\La_2^{(i)})-2l+1+s\,\ne\,0$
\,for all \,$s=1,\dots,l$, cf.~\cite{MV1}, \cite{MV2}, \cite{MTV3}.

\medskip
In the following, we need the next lemma.
\subsubsection{}
\vsk-.5>
\label{1221}
{\bf Lemma.}\enspace
{\it Let \,$m$ be a complex number and\/ \,$l$ a nonnegative integer.
\vvn.1>
Let \,$V$ \>be a $\glt$-module with weight decomposition
\vvn.1>
$V=\bigoplus_{k=0}^\infty V[k]$\,, \,where \>$V[k]\subset V$ \>is a weight
subspace of weight $(m-k,k)$. Assume that \,$m-2l+1+s\,\ne\,0$ \,for all
\,$s=1,\dots,l$. Then the map \,$e_{12}e_{21}:V[l-1]\to V[l-1]$ \,is
an isomorphism of vector spaces.}

\begin{proof}
Let \>$U_k=\ker\bigl(e_{12}^{l-k}|_{V[l-1]}^{\vp1}\bigr)$. Clearly,
\,$V[l-1]\,=\,U_0\supset U_1\dots\supset U_{l-1}\supset U_l\>=\>\{0\}$\;.

Let \,$C=e_{11}(e_{22}+1)-e_{12}e_{21}$.
Set \,$P(x)=\prod_{k=0}^{l-1}(x-c_k)$, \,where \,$c_k=k\>(m-k+1)$, \,and
\be
Q(x)\,=\;\frac{P(x)-P(c_l)}{x-c_l}\;.
\vv-.1>
\ee
We have
\vvn-.4>
\be
e_{12}e_{21}|_{V[l-1]}^{\vp1}\,=\;(c_l-C)|_{V[l-1]}^{\vp1}\;.
\ee

\smallskip
Since $C$ is a central element, we have \,$(C-c_k)\>U_k\subset U_{k+1}$.
\,Therefore, \,$P(C)|_{V[l-1]}=\>0$, \,and
\vvn-.3>
\be
(c_l-C)|_{V[l-1]}^{\vp1}\,Q(C)|_{V[l-1]}^{\vp1}\,=\,P(c_l)\;.
\vv.4>
\ee
The assumption on $m$ and $l$ implies that $P(c_l)\ne 0$.
\,Hence, the operator \,$(c_l-C)|_{V[l-1]}^{\vp1}$ is invertible.
\end{proof}

\subsection{Yangian}
\label{Yangian}
The Yangian $Y(\glt)$ is the unital associative algebra with
generators $T^{\{s\}}_{ab}$, $a,b=1,2$ and $s=1,2,\dots$. Let
\vvn-.3>
\be
T_{ab}(u)\ =\ \delta_{ab} + \sum_{s=1}^\infty\, T^{\{s\}}_{ab} u^{-s}\ ,
\qquad
a,b=1,2\ .
\ee
Then the defining relations in $\Y$ have the form
\vvn.2>
\beq
\label{defining rel}
(u-v)\,\bigl(T_{ab}(u)T_{cd}(v)-T_{cd}(v)T_{ab}(u)\bigr)\,=\,
T_{cb}(v)T_{ad}(u)-T_{cb}(u)T_{ad}(v)\ ,
\vv.2>
\eeq
for all $a,b,c,d$.
The Yangian is a Hopf algebra with coproduct
\beq
\label{Delta}
\Delta\ :\ T_{ab}(u) \ \mapsto \
\sum_{c=1}^2\, T_{cb}(u)\otimes T_{ac}(u)\
\vv-.6>
\eeq
for all $a,b$.

\subsubsection{}
\label{comm rel}
{\bf Proposition}\enspace
\cite{KBI}. \enspace
{\it The following relations hold\/}:
\begin{align*}
T_{11}(u)\,T_{12}(u_1)\dots T_{12}(u_k)\ &{}=\ \prod_{i=1}^k\,
\frac{u-u_i-1}{u-u_i}\ T_{12}(u_1)\dots T_{12}(u_k)\,T_{11}(u)\,+{}
\\
{}+\,\frac1{(k-1)!}\;T_{12}(u)\sum_{\si\in S_k\!}
& \,\biggl(\frac1{u-u_{\si_1}}\,\prod_{i=2}^k\,
\frac{u_{\si_1}-u_{\si_i}-1}{u_{\si_1}-u_{\si_i}}
\ T_{12}(u_{\si_2})\dots T_{12}(u_{\si_k})\,T_{11}(u_{\si_1})\biggr)\ ,
\end{align*}
\vvn-.8>
\begin{align*}
T_{22}(u)\,T_{12}(u_1)\dots T_{12}(u_k)\ &{}=\ \prod_{i=1}^k\,
\frac{u-u_i+1}{u-u_i}\ T_{12}(u_1)\dots T_{12}(u_k)\,T_{22}(u)\,+{}
\\
{}+\,\frac1{(k-1)!}\;T_{12}(u)\sum_{\si\in S_k\!}
& \,\biggl(\frac1{u-u_{\si_1}}\,\prod_{i=2}^k\,
\frac{u_{\si_1}-u_{\si_i}+1}{u_{\si_1}-u_{\si_i}}
\ T_{12}(u_{\si_2})\dots T_{12}(u_{\si_k})\,T_{22}(u_{\si_1})\biggr)\ .
\notag
\end{align*}

\subsubsection{}
\label{center}
A series $f(u)$ in $u^{-1}$ is called monic if \,$f(u)=1+O(u^{-1})$.
For a monic series $f(u)$, there is an automorphism
\be
\varphi_f\ :\ \Y\ \to\ \Y\ ,\qquad T_{ab}(u)\ \mapsto\ f(u)\,T_{ab}(u)\ .
\ee
There is a one-parameter family of automorphisms
\be
\rho_z \ :\ \Y\,\to\,\Y \qquad T_{ab}(u)\,\mapsto\,T_{ab}(u-z)\ ,
\ee
where in the right-hand side,
$(u-z)^{-1}$ has to be expanded as a power series in $u^{-1}$.

\smallskip
The Yangian $\Y$
contains the universal enveloping algebra $U(\glt)$ as a Hopf subalgebra.
The embedding is given by the formula $e_{ab} \mapsto T^{\{1\}}_{ba}$
for all $a,b$. We identify $U(\glt)$ with its image.

\smallskip
The evaluation homomorphism $\ep:\Y\to U(\glt)$ is defined by the rule:
$T^{\{1\}}_{ab} \mapsto e_{ba}$ \,for all $a,b$, \;and
\;$T^{\{s\}}_{ab}\mapsto 0$ \;for all $a,b$ and all $s>1$.

\smallskip
We denote by $^+ : \Y \to \Y$ the antiinvolution defined by
\beq
\label{anti}
\bigl(T_{ab}(u)\bigr)^+\,=\,T_{ba}(u)\,.
\eeq

\subsection{Bethe subalgebra}
The series
\vvn.2>
\beq
\label{qdet}
\qdet\,T(u)\ =
\ T_{1\,1}(u)\,T_{2\,2}(u-1)\,-\,T_{1\,2}(u)\,T_{2\,1}(u-1)\,
\vv.2>
\eeq
is called {\it the quantum determinant\/}. The coefficients of the series
\;$\qdet\,T(u)$ belong to the center of the Yangian $\Y$ \cite{IK}.

\medskip
The series $T_{11}(u)+T_{22}(u)$ is called the {\it transfer matrix\/}.
It is known that the coefficients of the series $T_{11}(u)+T_{22}(u)$
commute~\cite{FT}.

\medskip
We call the unital subalgebra $\Bg\subset\Y$ generated by coefficients
of the series $\qdet\,T(u)$ and $T_{11}(u)+T_{22}(u)$ {\it the Bethe
subalgebra\/}. The Bethe subalgebra is commutative. Elements of the Bethe
subalgebra commute with elements of the subalgebra $U(\glt)$ and are invariant
under the antiinvolution~\Ref{anti}.

\goodbreak
\subsection{Yangian modules}
\subsubsection{}
\vsk-.4>
{\bf Theorem}\enspace
\cite{T}.\enspace
{\it Let $V$ be an irreducible finite-dimensional $\Y$-module.
There exists a unique up to proportionality vector $v\in V$,
monic series $c_1(u),\>c_2(u)$, and a monic polynomial $\Po(u)$ such that
\begin{align*}
T_{21}(u)\, v\ &{}=\ 0 \ ,
\\
T_{aa}(u)\, v\ &{}=\ c_a(u)\, v \ , \qquad a=1, 2\ ,
\\[-20pt]
\end{align*}
and}
\vvn-.6>
\beq
\label{ccP}
\frac{c_1(u)}{c_{2}(u)}\;=\;\frac{\Po(u+1)}{\Po(u)}\;.
\eeq

\medskip
The vector $v$ is called {\it a highest weight vector\/}, the series
\>$c_1(u)\>,\,c_2(u)$ --- {\it the Yangian highest weights\/}, and
the polynomial $\Po(u)$ --- {\it the Drinfeld polynomial\/} of the module $V$.

\subsubsection{}
{\bf Theorem}\enspace
\cite{T}.\enspace
{\it For any monic series $c_1(u),\>c_2(u)$ and a monic polynomial $\Po(u)$
obeying relation \Ref{ccP}, there exists a unique irreducible
finite-dimensional $\Y$-module $V$ such that \,$c_1(u)\>,\,c_2(u)$ are
the Yangian highest weights of the module $V$.}

\subsubsection{}
Let $V_1,V_2$ be irreducible finite-dimensional $\Y$-modules with respective
highest weight vectors $v_1,v_2$. Then for the $\Y$-module $V_1\otimes V_2$,
we have
\begin{align*}
T_{21}(u)\, v_1\otimes v_2\ &{}=\ 0 \ ,
\\
T_{aa}(u)\, v_1\otimes v_2\ &{}=
\ c_a^{(1)}(u)\,c_a^{(2)}(u)\, v_1\otimes v_2\ , \qquad a=1,2\ .
\end{align*}

Let $W$ be the irreducible subquotient of $V_1\otimes V_2$ generated by
the vector $v_1\otimes v_2$. Then the Drinfeld polynomial of the module $W$
equals the products of the Drinfeld polynomials of the modules
$V_1$ and $V_2$.

\subsubsection{}
For a $\glt$-module $V$, let the $\Y$-module $V(z)$ be the pullback of $V$
through the homomorphism $\ep\circ\rho_z$\>; that is, the series $T_{ab}(u)$
acts on $V(z)$ as $1+(u-z)^{-1}e_{ba}$\,. The module $V(z)$ is called
the evaluation module with evaluation point $z$.

\subsubsection{}
Let $\bs \La = (\La^{(1)},\dots,\La^{(n)})$ be a collection of integral
dominant $\glt$-weights, where $\La^{(i)}=(\La_1^{(i)},\La_2^{(i)})$,
\,$\La_1^{(i)}\ge\La_2^{(i)}$, for $i=1,\dots, n$. For generic complex numbers
$z_1,\dots,z_n$, the tensor product of evaluation modules
\be
L_{\bs \La}(\bs z)\ =
\ L_{\La^{(1)}}(z_1) \otimes\dots\otimes L_{\La^{(n)}}(z_n)
\vv.3>
\ee
is an irreducible finite-dimensional $\Y$-module and the corresponding highest
weight series $c_1(u)\>,\,c_2(u)$ have the form
\vvn-.4>
\beq
\label{cau}
c_a(u) \ =\ \prod_{i=1}^n\,\frac {u-z_i + \La_a^{(i)}}
{u-z_i }\ .
\vv.2>
\eeq
The corresponding Drinfeld polynomial equals
\vvn.2>
\be
\Po(u)\,=\,\prod_{i=1}^n\;
\prod_{s=\La_{2}^{(i)}}^{\La_1^{(i)}-1}(u-z_i+s)\, .
\ee

\subsection{Holomorphic representation}
\label{holomorphic}
The results of this section go back to~\cite{T}.

\smallskip
Choose monic polynomials $a(u)\>,\,d(u) \in \C[u]$ of positive degree $n$,
\beq
\label{adu}
a(u)\,=\,\prod_{i=1}^n\,(u-z_i+m_i)\,, \qquad
d(u)\,=\,\prod_{i=1}^n\,(u-z_i)\,.
\eeq

\subsubsection{}
{\bf Proposition.}
{\it There exists a unique $\Y$-action on the vector space $\Cx$ such that
\vvn-.2>
\begin{align}
\label{T12}
\bigl(T_{12}(u)\cdot p\bigr)(x)\ & {}=
\ \frac1{d(u)}\,\sum_{i=1}^n\, u^{n-i}\, x_i\,p(x_1,\dots,x_n)
\\[4pt]
&{}=\ \biggl(\frac {x_1}u +\frac {x_2-x_1\sum_{i=1}^nz_i}{u^2} + \dots{}\biggr)
\,p(x_1,\dots,x_n)
\notag
\end{align}
for any polynomial \;$p\in\C[x_1,\dots,x_n]$\,, \;and
\vvn.3>
\beq
\label{Tab1}
T_{11}(u)\,\cdot\,1\ =\ \frac{a(u)}{d(u)}\,\cdot\,1\ ,
\qquad T_{22}(u)\,\cdot\,1\ =\ 1\ ,\qquad T_{21}(u)\,\cdot\,1\ =\ 0\ ,
\vv.2>
\eeq
where $1$ stands for the constant polynomial equal to $1$ as an element of
$\Cx$.}

\medskip
We denote by $\W$\, the $\Y$-module defined by formulae~\Ref{T12}, \Ref{Tab1}
and call it {\it the holomorphic representation\/} of $\Y$,
associated with the polynomials $a(u)\>,\,d(u)$.

\medskip
The Yangian module $\W$ is cyclic: every element of $\W$ can be obtained from
$1$ by the action of a suitable polynomial in $T^{\{1\}}_{12},$
$T^{\{2\}}_{12},\;\dots{}\;$. Formulae~\Ref{Tab1} mean that $1$ is
an eigenvector of the operators $T^{\{s\}}_{11},T^{\{s\}}_{22}$ and $1$ is
annihilated by the operators $T^{\{s\}}_{21}$ with $s=1,2,\dots{}$\,.
Then the Yangian commutation relations \Ref{defining rel} allow us to determine
the action of $T^{\{s\}}_{11},T^{\{s\}}_{22}$, $T^{\{s\}}_{21}$ on all elements
of $\W$.

\smallskip
Since the coefficients of the series \,$\qdet T(u)$ are central,
and the module $\W$ is generated by the polynomial $1$, we have
\beq
\label{qdetsing}
\qdet\,T(u)\big|_{\W}\,=\;\frac{a(u)}{d(u)}\;.
\eeq

For every $i,j = 1,2$, we have
\vvn-.6>
\beq
\label{def of tilde-T}
T_{ij}(u)\big|_{\W}\,=\;\frac{\tT_{ij}(u)}{d(u)}\;,
\eeq
where $\tT_{ij}(u)$ is an $\End\,(\W)$-valued polynomial in $u$ of degree \,$n$
\>for \,$i=j$, and of degree \,$n-1$ \,for \,$i\ne j$.

\subsubsection{}
The embedding $U(\glt) \hookrightarrow \Y$ defines a $\glt$-module structure
on $\W$. The $\glt$-weight decomposition of $\W$ is the degree decomposition
$\W\,=\,\oplus_{l=0}^\infty\,\W[\>l\>]$ into subspaces of homogeneous
polynomials. The subspace $\W[\>l\>]$ of homogeneous polynomials of degree $l$
has $\glt$-weight $\bigl(\>\sum_{i=1}^n m_i-l,\,l\>\bigr)$.

\subsubsection{}
\label{singular}
{\bf Lemma.}
{\it Let \ ${\Sing\,\W[\>l\>]\,=\,\{\,p\in\W[\>l\>]\ |\ e_{12} p\,=\,0\,\}}$
\;be the subspace of \,$\glt$-sin\-ular vectors. Assume that the pair
$((m_1,0),\dots,(m_n,0))$, $l$ is separating. Then}
\vvn.2>
\be
\dim\,\Sing\,\W[\>l\>]\ =\ \dim\,\W[\>l\>]\ -\ \dim\,\W[\>l-1\>]\ .
\ee

\begin{proof}
The map \,$e_{12}e_{21} : \W[\>l-1\>] \to \W[\>l-1\>]$
\,is an isomorphism of vector spaces since the pair
$((m_1,0),\dots,(m_n,0))$, $l$ is separating, see Lemma~\ref{1221}.
The fact that $e_{12}e_{21}$ is an isomorphism implies the lemma.
\end{proof}

\subsubsection{}
Denote by ${}^+ : \Y \to \Y$ the antiinvolution defined by
$T_{ij}^+(u) = T_{ji}(u)$.
Denote by $\,\phi : \W \to \C\,$
the linear function \, $p(x_1,\dots,x_n) \mapsto p(0,\dots,0)$.
{\it The Yangian Shapovalov form\/} on $\W$ is the unique symmetric bilinear form $S$
on $\W$ defined by the formula $S(x\cdot 1,\,y\cdot 1)\ =\ \phi ( x^+y\cdot 1)$
\,for all $x,y \in \Y$.

\medskip
Different $\glt$-weight subspaces of $\W$ are orthogonal with respect to
the form $S$, and
\be
\det\,S|_{\W[\>l\>]}\ =\ {\rm const}\, \prod_{i,j = 1}^n \prod_{s=0}^{l-1}\,
(z_i-z_j+m_j-s){\vphantom{\big|}}^{\binom{n+l-s-2}{n-1}}\,,
\ee
where the constant does not depend on
$z_1,\dots,z_n$\>, \>$m_1,\dots,m_n$.

\subsection{}
The kernel of the Yangian Shapovalov form $K\subset \W$ is a $\Y$-submodule.
The $\Y$-module $\W/K$ is irreducible.
The Yangian Shapovalov form on $\W$ induces a nondegenerate symmetric bilinear
form on $\W/K$ called the {\it Yangian Shapovalov form} of the module $\W/K$.

\subsubsection{}
\label{sec product of Verma}
{\bf Theorem} \cite{T}.
{\it For generic $z_1,\dots,z_n$\>, \>$m_1,\dots,m_n$, the $\Y$-module $\W$
is irreducible and isomorphic to the tensor product of evaluation Verma modules
$M_{(m_1,0)}(z_1)\otimes \dots \otimes M_{(m_n,0)}(z_n)$.
Any such an isomorphism sends $1$ to a scalar multiple of the tensor product
$v_{(m_1,0)}\otimes \dots \otimes v_{(m_n,0)}$ of highest weight vectors of
the corresponding Verma modules.}

\subsubsection{}
\label{sec product of irrep}
{\bf Theorem} \cite{T}.
{\it Let \,${m_i\in\Z_{\geq 0}}$ for \,${i=1,\dots,n}$, \>and
\,${m_1\leq m_2\leq\dots\leq m_n}$. Assume that \,$z_i-z_j+m_j-s\ne 0$
\,and \,$z_i-z_j-1-s\ne 0$ \,for all \,$i<j$ and \,$s=0,1,\dots,\alb m_i-1$.
Then for any permutation $\si\in S_n$, the
irreducible $\Y$-module $\W/K$ is isomorphic to the tensor product of
evaluation irreducible modules
$L_{(m_{\si_1},0)}(z_{\si_1})\otimes\dots\otimes
L_{(m_{\si_n},0)}(z_{\si_n})$. Any such an isomorphism sends the element
corresponding to $1$ to a scalar multiple of the tensor product
$v_{(m_{\si_1},0)}\otimes\dots\otimes v_{(m_{\si_n},0)}$ of highest
weight vectors of the corresponding irreducible modules. }

\medskip
For a proof of this theorem see also \cite{CP}.

\subsubsection{}
The assumption of Theorem \ref{sec product of irrep} can be formulated
geometrically as the assumption that for $i< j$ the sets
\;$Z_i=\{ z_i, \,z_i-1, \,\dots\,,z_i-m_i\}$ \;and
\;$Z_j\ =\ \{ z_j,\, z_j-1,\, \dots\,,\alb z_j-m_j\}$
\;either do not intersect, or the smaller set $Z_i$ is a subset
of the larger set $Z_j$ (since we assumed that $m_i\le m_j$).

\section{Algebra $A_W$ and universal difference operator}
\label{Universal difference operator}

\subsection{Definition}
Let $V$ be a $\Y$-module. We call the image of the Bethe algebra $\Bg\subset\Y$
in $\End\,(V)$ the {\it Bethe algebra\/} associated with $V$. If $U\subset V$
is a vector subspace preserved by elements of the Bethe algebra $\Bg_V$, then
their restrictions to $U$ define a commutative unital subalgebra
$\Bg_U\subset\End\,(U)$ called {\it the Bethe algebra\/} associated with $U$.

\subsubsection{}
Define the operator $\shift$ acting on functions of $u$ as
$(\shift f)(u)=f(u+1)$.

\smallskip
Let $V$ be a $\Y$-module such that for all $a,b$ the series $T_{ab}(u)|_V$
sum up to $\End\,(V)$-valued rational functions in $u$. Let $U\subset V$ be
a vector subspace preserved by the Bethe algebra $\Bg_V$. The {\it universal
difference operator\/} $\DD_U$ acting on $U$-valued functions in $u$ is defined
by the formula
\vvn.3>
\be
\DD_U\;=\;1\,-\,\bigl(T_{11}(u)+T_{22}(u)\bigr)\big|_U\,\shift^{-1}\,+\,
\qdet\,T(u)\big|_U\,\shift^{-2}\;,
\vv.3>
\ee
see \cite{Tal}, \cite[\!(4.16)\,]{MTV1}, \cite{MTV2}.
The operator $\DD_U$ is a linear second-order difference operator.

\subsection{Algebra $A_W$}{Operator $\DD_{\>\Sing\,\W[\>l\>]}$}
\label{sec nash M}
Consider the Bethe algebra $\Bg_{W_{a,d}}$ associated with the $\Y$-module
$W_{a,d}$. Recall that
\be
\bigl(\qdet\,T(u)\bigr)\big|_{\,\W} \,=\;\frac{a(u)}{d(u)}\;,
\ee
see~\Ref{qdetsing}, and
\vvn-.5>
\be
\bigl(T_{11}(u) + T_{22}(u)\bigr)\big|_{\,\W}\ =
\ \frac{B(u,\bs{\tilde H})}{d(u)}
\ee
where
\vvn-.5>
\beq
B(u,\bs{\tilde H})\,=\,\tilde H_0u^n + \tilde H_1u^{n-1} + \dots + \tilde H_n
\vvn.5>
\eeq
for suitable coefficients $\tilde H_k \in \End\,\bigl(\W\bigr)$,
see~\ref{def of tilde-T}. It follows from Proposition~\Ref{comm rel} that
the coefficients $\tilde H_0\>,\,\tilde H_1$ are scalar operators,
\,$\tilde H_0=2$\>, \ $\tilde H_1 = \sum_{i=1}^n (m_i- 2 z_i)$.

The elements $\tilde H_k$ are called {\it the \XXX/ Hamiltonians}
associated with $\W$.

\subsubsection{}
\label{H-Hamiltonians}
The Hamiltonians $\tilde H_k$ preserve the subspace $\Sing\,\W[\>l\>]$
defined in Section \ref{singular}. Set
\vvn-.1>
\be
H_k\,=\,\tilde H_k|_{\,\Sing\,\W[\>l\>]}\in\End\,(\Sing\,\W[\>l\>])
\vv-.4>
\ee
and
\be
B(u, \bs H)\,=\,H_0u^n + H_1u^{n-1} + \dots + H_n\;.
\vvgood
\vv.4>
\ee
The coefficients $H_0\>,\,H_1\>,\,H_2$ are scalar operators,
\begin{align*}
& H_0\ =\ 2\ , \quad\qquad H_1\ =\ \sum_{i=1}^n\, (m_i- 2 z_i)\ ,
\\
& H_2\ =\ l\,\bigl(l - 1 -\sum_{i=1}^n\,m_i\bigr)\ +\
\sum_{1\leq i< j\leq n}\!\bigl(z_iz_j + (z_i-m_i)(z_j-m_j)\bigr)\ .
\end{align*}
The simplest way to get the last formula is to extract $H_2$ from
the coefficient of $u^{-2}$ of the series
\,$\qdet\,T(u)\big|_{\,\Sing\,\W[\>l\>]}$\>, \>see~\Ref{qdet}, and
to use formula~\Ref{qdetsing}.

\smallskip
We denote by $A_W$ the Bethe algebra associated with $\Sing\,\W[\>l\>]$.
It is the unital subalgebra of $\End\,\bigl(\Sing\,\W[\>l\>]\bigr)$
generated by the operators $ H_3, H_4, \dots , H_n$, called {\it the
\XXX/ Hamiltonians} associated with $\Sing\,\W[\>l\>]$.

\subsubsection{}
\label{Algebra A_W and Shapovalov form}
The operators of the algebra $A_W$ are symmetric with respect
to the Yangian Shapovalov form on $\W$,
\vvn-.2>
\be
S(fv,w)\,=\,S(v, fw)
\vv.3>
\ee
for all $f\in A_W$ and $v,w \in \W$, see \cite{MTV1}.

\subsection{Operator $\D_{\>\Sing\,\W[\>l\>]}$}
Consider the universal difference operator $\DD_{\>\Sing\,\W[\>l\>]}$
acting on $\Sing\,\W[\>l\>]$-valued functions,
\vvn-.2>
\be
\DD_{\>\Sing\,\W[\>l\>]}\,=\;1\,-\,\frac {B(u, \bs H)}{d(u)}\,\shift^{-1}
\,+\,\frac {a(u)}{d(u)} \, \shift^{-2}\,,
\vv.2>
\ee
The {\it modified universal difference operator\/}
$\D_{\>\Sing\,\W[\>l\>]}$ is defined by the formula
\vvn.4>
\be
\D_{\>\Sing\,\W[\>l\>]}\,=\,d(u)\,\DD_{\>\Sing\,\W[\>l\>]}\,.
\vvn-.6>
\ee
Then
\vv-.2>
\be
\D_{\>\Sing\,\W[\>l\>]}\;=
\;d(u)\>-\>B(u, \bs H)\,\shift^{-1}\>+\>a(u) \, \shift^{-2}\ .
\ee

\subsubsection{}
\label{thm exist of polyn}
{\bf Theorem.}
{\it
Assume that the pair $((m_1,0),\dots, (m_n,0)), \, l$\ is separating.
Then for any $v_0\in \snash$ there exist unique $v_1,\dots,v_l \in\snash$
such that the function
\be
w(u) \ =\ v_0\,u^l \, +\,v_1\,u^{l-1}\, +\, \dots \,+\,v_l
\ee
is a solution of the difference equation $\D_{\>\snash} w(u)\,=\,0$.}

\begin{proof}
By Lemma~\ref{singular} the dimension of $\snash$ does not depend on
$z_1,\dots,z_n$\>, \>$m_1,\dots,m_n$, if the pair
$\bigl((m_1,0),\dots,(m_n,0)\bigr)\>,\,l$ is separating. Because of that,
we may consider the difference equation $\D_{\>\Sing\,\W[\>l\>]} v(u)\,=\,0$
as a difference equation on a fixed vector space with coefficients of
the difference equation algebraically depending on parameters
$z_1,\dots,z_n$\>, \>$m_1,\dots,m_n$.

Given a vector $v_0\in\snash$, we look for a solution of the difference
equation $\D_{\>\Sing\,\W[\>l\>]} v(u)\,=\,0$ in the form
$v_0u^l + \sum_{j=1}^\infty v_ju^{l-j}$. Substituting this expression into
the equation, we can calculate all of the coefficients $v_j$ recursively,
and they are algebraic functions of $z_1,\dots,z_n$\>, \>$m_1,\dots,m_n$.

For generic $z_1,\dots,z_n$ \>and large positive integral \>$m_1,\dots,m_n$\>,
the coefficients $v_j$ are equal to zero for all $j>l$ by
Theorem~\ref{sec product of irrep} and~\cite[Theorem~7.3]{MTV3}.
Hence, the same coefficients are equal to zero for all $z_1,\dots,z_n$\>,
\>$m_1,\dots,m_n$ such that the pair $\bigl((m_1,0),\dots,(m_n,0)\bigr)\>,\,l$
\,is separating.
\end{proof}

\section{Algebra $A_D$}
\label{Algebra A_D}

\subsection{Definition}
\label{subsub h's}
{}From now on until the end of Section~\ref{Equations with polynomial solutions
only} we fix complex numbers $z_1,\dots,z_n$\>, $m_1,\dots,m_n$\>, and
a nonnegative integer \>$l$. We always assume that the polynomials $a(u)$
and $d(u)$ are given by formulae~\Ref{adu}.

\medskip
Let $\bs a = (a_1,\dots,a_l)$ and $\bs h = (h_1,\dots,h_n)$.
Consider the space $\C^{l+n}$ with coordinates $\bs a, \bs h$.
Let $D$ be the affine subspace of $\C^{l+n}$ defined by equations
$q_{1}(\bs h)=0$\>, \>$q_2(\bs h) = 0$, where
\vvn-.3>
\begin{align*}
q_{1}(\bs h)\,&{}=\,h_1\>-\>\sum^n _{i=1}\, (m_i-2z_i) \;,
\\
q_2(\bs h)\,&{}=\,h_2\>-\,
l\,(l - 1 -\sum_{i=1}^n\,m_i)\,-
\sum_{1\leq i< j\leq n}\,(z_iz_j + (z_i-m_i)(z_j-m_j))\;.
\end{align*}
Let
\vvn-.8>
\begin{gather}
\label{D}
\qa\,=\,u^l + a_1u^{l-1}+ \dots + a_l\;,
\\[4pt]
B(u, \bs h)\,=\,2u^n + h_1u^{n-1} + \dots + h_n \;,
\notag
\\[4pt]
\Dh\,=\,d(u)\>-\>B(u, \bs h)\,\shift^{-1}\>+\>a(u)\,\shift^{-2}\;.
\notag
\end{gather}
If $\bs h$ satisfy the equations $q_{1}(\bs h) = 0$ and $q_{2}(\bs h) = 0$,
then the polynomial $\Dh(\qa)$ is a polynomial in $u$ of degree $l+n-3$,
\vvn.3>
\be
\Dh(\qa)\,=\,q_3(\bs a,\bs h)\,u^{l+n-3}\>+\ \dots \ +\>q_{l+n}(\bs a,\bs h)\;.
\vv.3>
\ee
The coefficients $q_i (\bs a,\bs h)$ are linear functions in $\bs a$ and
linear functions in $\bs h$.

Denote by $I_D$ the ideal in
$\C[\bs a, \bs h]$ generated by
polynomials $q_{1}, q_2, q_3,\dots,q_{l+n}$.
The ideal $I_D$ defines a scheme $C_D\subset D$. Then
\bea
A_D\ =\ \C[\bs a, \bs h]/ I_D\
\eea
is the algebra of functions on $C_D$. The scheme $C_D$ is the scheme of points
$\T\in D$ such that the polynomial $p(u,\bs a(\T))$ solves the difference
equation $\Dhp w(u)=0$.

\subsection{Independence of dimension of $A_D$ on $z_1,\dots, z_n$}

For fixed $m_1,\dots,m_n$, the scheme $C_D$ and the algebra $A_D$
depend on the choice of numbers $\bs z = (z_1,\dots,z_n)$:
$C_D = C_D(\bs z)$,\ $A_D = A_D(\bs z)$.

\subsubsection{}
\label{thm const wrt z}
{\bf Theorem.}
{\it
Assume that the pair $\bigl((m_1,0),\dots, (m_n,0)\bigr)\>,\,l$ is separating.
Then the dimension of $A_D(\bs z)$, considered as a vector space, is finite
and does not depend on the choice of numbers $z_1,\dots,z_n$.}

\goodbreak
\begin{proof}
It suffices to prove two facts:
\begin{enumerate}
\item[(i)]
For any $\bs z$, there are no algebraic curves over $\C$ lying in $C_D(\bs z)$.

\item[(ii)]
Let a sequence $\bs z^{(i)}$, $i=1,2,\dots{}$,
tend to a finite limit $\bs z = (z_1,\dots,z_n)$.
Let $\T^{(i)} \in C_D(\bs z^{(i)})$\>, \,$i = 1, 2, \dots\ $,
be a sequence of points. Then all coordinates
$\bigl(\bs a(\T^{(i)}), \bs h(\T^{(i)}\bigr)$
remain bounded as $i$ tends to infinity.
\end{enumerate}
By fact~(i), the dimension of $A_D(\bs z)$ is finite for any $\bs z$, whereas
fact~(ii) implies that \,$\dim A_D(\bs z)$ does not depend on $z_1,\dots,z_n$.

\medskip
For a point $\T$ in $C_D(\bs z)$, the operator $\D_{\bs h(\T)}$ has the form
\be
d(u)\ - \ (2u^n + h_1(\T)u^{n-1} + h_2(\T)u^{n-2} + h_3(\T)u^{n-3}+
\dots + h_n(\T))\,\shift^{-1}\ + \ a(u)\,\shift^{-2}\ ,
\ee
where the coefficients $h_1(\T)\>,\,h_2(\T)$ are determined by the equations
\>$q_{1}(\bs h) = 0$ \>and \>$q_{2}(\bs h) = 0$.

Assume that (i) is not true. Since any affine algebraic curve over $\C$
is unbounded,
there exists a sequence of points $\T^{(i)} \in C_D(\bs z)$, $i=1,2,\dots{}$,
which tends to infinity as $i$ tends to infinity. Then it is easy to see that
$\bs h(\T^{(i)})$ cannot tend to infinity since it would contradict the fact
that $\D_{\bs h(\T^{(i)})}\bigl(p(u,\bs a(\T^{(i)}))\bigr) = 0$.
Choosing a subsequence, we may assume that $\bs h(\T^{(i)})$ has a finite limit
as $i$ tends to infinity. Then $\bs a(\T^{(i)})$ cannot tend to infinity since
it would mean that the limiting difference equation has a polynomial solution
of degree less than $l$, and this is impossible.

This reasoning implies that $\T^{(i)}\in C_D(\bs z)$ cannot tend to infinity.
Thus we get a contradiction and statement~(i) is proved.

\smallskip
The proof of statement~(ii) is similar.
\end{proof}

\subsection{Second description of $A_D$ and
epimorphism $\psi_{DW} : A_D \to A_W$}
\label{Epimorphism psi-DM}

\subsubsection{}
\vsk-.5>
\label{thm-conjecture 2}
{\bf Theorem.}
{\it
Assume that the pair $\bigl((m_1,0),\dots, (m_n,0)\bigr)\>,\,l$ is separating.
Assume that $\bs h$ satisfies equations $q_{1}(\bs h) = 0$ and
$q_{2}(\bs h) = 0$. Consider the system
\beq
\label{system q}
q_i(\bs a, \bs h)\ =\ 0\ , \qquad i = 3, \dots , l+2\ ,
\eeq
as a system of linear equations with respect to $a_1,\dots,a_l$.
Then this system has a unique solution \ $a_i=a_i(\bs h)\>,\ \,i=1,\dots,l$,
\,where $a_i(\bs h)$ are polynomials in \,$\bs h$.}
\qed

\begin{proof}
The claim follows from the fact that
\be
q_{2+i}(\bs a,\bs h)\,=\,
i\,\Bigl(\>\sum_{s=1}^n m_s - 2l + i + 1 \Bigr)\, a_i \,+
\>\sum_{j=1}^{i-1} \,q_{ij}(\bs h)\, a_j
\ee
for $i = 1, \dots , l$. Here $q_{ij}$ are some linear functions of $\bs h$.
The coefficient of $a_i$ does not vanish since the pair
$\bigl((m_1,0),\dots, (m_n,0)\bigr)\>,\,l$ is separating.
\end{proof}

\subsubsection{}
Denote by $I'_D$ be the ideal in $\C[\bs h]$ generated by the polynomials
$q_{1}(\bs h)\>,\,q_2(\bs h)$\>,\alb\ $q_j(\bs a(\bs h),\bs h)$\>,
\,$j= l+3,\dots, l+n$. The ideal $I'_D$ defines a scheme $C'_D$ in the space
$\C^n$ with coordinates $\bs h=(h_1,\dots,h_n)$. The scheme $C'_D$ is
the scheme of points $\bs r\in\C^n$ such that the difference equation
$\D_{\bs h(\bs r)}w(u)=0$ has a polynomial solution of degree $l$.

\smallskip
Theorem~\ref{thm-conjecture 2} implies that
\vvn-.2>
\beq
\label{sec 2nd description}
A_D \ \cong \ \C[\bs h]/ I'_D\ .
\vv.2>
\eeq

\smallskip
Let $H_{1},\dots,H_n$ be the operators introduced in
Section~\ref{H-Hamiltonians}.

\subsubsection{}
\label{thm D to B}
{\bf Theorem.}
{\it
Assume that the pair $\bigl((m_1,0),\dots, (m_n,0)\bigr)\>,\,l$ is separating.
Then the assignment $h_s\ \mapsto\ H_s$\>, \;$s = 1, \dots , n$, determines
an algebra epimorphism \,$\psi_{DW} : A_D \to A_W$.
}

\begin{proof}
We use description~\Ref{sec 2nd description} of the algebra $A_D$.
The equations defining the scheme \,$C'_D$ are the equations of existence
of a polynomial solution of degree $l$ to the polynomial difference equation
$\D_{\bs h}w(u)=0$. The operators $H_{1},\dots,H_n$ satisfy the defining
equations for \,$C'_D$ by Theorem~\ref{thm exist of polyn}.
\end{proof}

\section{Bethe ansatz equations}
\label{Bethe ansatz and C_D}

\subsection{Bethe ansatz equations}
The Bethe ansatz equations is the following system of equations
with respect to complex numbers $\bs t=(t_1,\dots,t_l)$\,:
\begin{gather}
\label{BAE}
\prod_{s=1}^n\,(t_j - z_s + 1 + m_s)\,\prod_{k \ne j}\,(t_j - t_k - 1)\;=\;
\prod_{s=1}^n\,(t_j-z_s +1)\,\prod_{k \ne j}\,(t_j - t_k + 1)\;,
\\
\rightline{$j = 1, \dots , l\;.$\phantom{(4.1)}}
\notag
\end{gather}
A solution $\bs t$ is called admissible if all $t_1,\dots,t_l$ are distinct,
and all factors in \Ref{BAE} are nonzero.

The permutation group $S_l$ acts on admissible solutions.
If $\bs t=(t_1,\dots,t_l)$ is an admissible solution,
then any permutation of these numbers is an admissible solution too.
We shall consider $S_l$-orbits of admissible solutions.

The following lemma is well-known, see for example Lemma 2.2 in \cite{MV2}.

\subsubsection{}
\label{lemma on bethe ansatz}
{\bf Lemma.}
{\it
Let $\bs t$ be an admissible solution of system~\Ref{BAE}. Let
\be
p(u)\ =\ \prod_{i=1}^l(u-t_i)\ ,
\qquad \Bc(u)\ =\ \frac{d(u) p(u) + a(u) p (u-2) }{p(u-1)}\ .
\vv.1>
\ee
Then $\Bc(u)$ is a polynomial of degree $n$ and $p(u)$ is annihilated
by the difference operator}
\vvn.3>
\be
d(u)\ - \ \Bc(u)\,\shift^{-1}\
+ \ a(u)\,\shift^{-2}\ .
\vv->
\ee
\qed

\subsubsection{}
\label{cor on bethe ansatz}
{\bf Corollary.}
{\it
Any $S_l$-orbit of admissible solutions of the Bethe ansatz equations gives
a point of the scheme $C_D$ considered as a set. Moreover, different
$S_l$-orbits give different points.}
\qed

\subsubsection{}
\label{thm BAE have many slns}
{\bf Theorem.}
{\it
Assume that
the pair $\bigl((m_1,0),\dots,\alb(m_n,0)\bigr)\>,\,l$ is separating.
Then for generic $z_1,\dots, z_n$ the Bethe ansatz equations
have at least \,$\dim\,\Sing\,\W[\>l\>]$ distinct $S_l$-orbits of admissible
solutions.}

\subsubsection{}
\label{cor C_D has many points}
{\bf Corollary.}
{\it
Assume that the pair $\bigl((m_1,0),\dots,(m_n,0)\bigr)\>,\,l$ is separating.
Then for generic $z_1,\dots, z_n$ the scheme $C_D$ considered as a set has
at least \,$\dim\,\Sing\,\W[\>l\>]$ distinct points.}
\qed

\begin{proof}[Proof of Theorem \ref{thm BAE have many slns}]
Make the change of variables: $z_s = \hat z_s/\eps$\>, \,$s = 1, \dots, n$, and
$t_i=\hat t_i/\eps$\>, \,$i=1,\dots, l$. Then equations \Ref{BAE} take the form
\vvn.3>
\beq
\label{BAE-1}
\prod_{s=1}^n\,
\frac{\hat t_j - \hat z_s + \eps + m_s\eps}{\hat t_j-\hat z_s +\eps}
\ \prod_{k \ne j}\,\frac{\hat t_j - \hat t_k -\eps}{\hat t_j - \hat t_k +\eps}
\ = \ 1\ , \qquad j = 1, \dots , l\ .
\vv.2>
\eeq
As $\eps$ tends to zero, equations \Ref{BAE-1} take the form
\vvn.3>
\be
\sum_{s=1}^n \frac {m_s} {\hat t_j-\hat z_s} \ -
\ \sum_{k \ne j} \frac {2}{\hat t_j - \hat t_k } \ = \ O(\eps) \ ,
\qquad j = 1, \dots , l\ ,
\vv.2>
\ee
and in the limit we obtain
\vvn.3>
\beq
\label{BAE-3}
\sum_{s=1}^n \frac {m_s}{\hat t_j-\hat z_s} \ -
\ \sum_{k \ne j}\frac {2}{\hat t_j - \hat t_k } \ = \ 0 \ ,
\qquad j = 1, \dots , l\ .
\vv.2>
\eeq

The last system is the system of the Bethe ansatz equations for
the Gaudin model. It was proved in \cite{RV} that if the pair
$\bigl((m_1,0),\dots,(m_n,0)\bigr)\>,\,l$ is separating and
$\hat z_1, \dots, \hat z_n$ are generic, then system~\Ref{BAE-3}
has at least \,$\dim\,\Sing\,\W[\>l\>]$ \>distinct $S_l$-orbits of admissible
solutions. This proves Theorem \ref{thm BAE have many slns}.
\end{proof}

\section{Separation of variables}
\label{On separation for slt}

\subsection{Change of variables}
\label{change of var slt}
For a nonnegative integer $l$ let
$\Cl$ be
the vector space of symmetric polynomials in $y_1,\dots,y_{n-1}$
of degree not greater than $l$ with respect to each variable.
Let
\vvn.3>
\be
\Wy[\>l\>]\, =\, y_0^l \,\Cl\ \subset\ \C[y_0,y_1,\dots,y_{n-1}]
\vv.3>
\ee
and set $\Wy = \oplus_{l=0}^\infty \Wy[\>l\>]$.
Define an isomorphism of vector spaces
\vvn.3>
\beq
\label{isom of separ}
\W\;\cong\;\Wy
\vv-.2>
\eeq
using the formula
\vvn-.4>
\be
\sum_{i=1}^n\,x_i u^{n-i}\ = \
y_0\,\prod_{j=1}^{n-1}\,(u-y_j)\ ,
\ee
that is, by setting
\be
x_i \ = \ (-1)^{i-1}\,y_0\, \sigma_{i-1}(y_1,\dots,y_{n-1})\ ,
\ee
where $\sigma_{i-1}$ is the $(i-1)$-st elementary symmetric function.
For example, for $n=2$ \,we have \,$x_1u+x_2 = y_0 (u-y_1)$ \,and
\,$x_1= y_0$\>, \,$x_2 = -y_0y_1$.

\medskip
We will identify the spaces $\W$ and $\Wy$ using isomorphism
\Ref{isom of separ}. In particular, this defines a $\Y$-module structure
on $\Wy$. We denote by $\Sing\,\Wy[\>l\>] \subset \Wy[\>l\>]$ the subspace of
$\glt$-singular vectors.

Isomorphism \Ref{isom of separ} defines on $\Wy$ and its subspaces
the operators which were previously defined on $\W$ and its subspaces.
Those new operators will be denoted by the same symbols. In particular,
we shall consider the action of operators $\tT_{ij}(u)$ and
$\tilde H_0, \dots, \tilde H_n$ on $\Wy$.

\subsection{Sklyanin's theorem}
\subsubsection{}
\vsk-.6>
\label{Sklyanin thm}
{\bf Theorem \cite{Sk}.}
{\it
The action of $e_{11}, \,e_{22},\, \tT_{11}(u),\,
\tT_{22}(u)$ on $\Wy$
is given by the following formulae:
\begin{align}
\label{eii}
e_{11}\ &{}=\ \sum_{i=1}^n m_i\,-\,y_0\frac{\partial}{\partial y_0}\ ,
\qquad e_{22}\ =\ \,y_0\frac {\partial}{\partial y_0}\ ,
\\
\label{tT11}
\tT_{11}(u)\ &{}=\ \Bigl(u+e_{11}-\sum_{i=1}^n z_i+\sum_{j=1}^{n-1} y_j\Bigr)
\prod_{j=1}^{n-1} (u-y_j)\ +\ \sum_{j=1}^{n-1} a(y_j) \prod_{j' \ne j}
\frac {u-y_{j'}}{y_j-y_{j'}}\ \shift_{y_j}^{-1}\ ,
\\[4pt]
\label{tT22}
\tT_{22}(u)\ &{}=\ \Bigl(u+e_{22}-\sum_{i=1}^n z_i+\sum_{j=1}^{n-1} y_j\Bigr)
\prod_{j=1}^{n-1} (u-y_j)\ +\ \sum_{j=1}^{n-1} d(y_j) \prod_{j' \ne j}
\frac {u-y_{j'}}{y_j-y_{j'}}\ \shift_{y_j}\ ,
\end{align}
where $\shift_{y_j}^{}:f(y_0,\dots,y_{n-1})\mapsto
f(y_0,\dots,y_j + 1, \dots, y_{n-1})$.
}
\begin{proof}
The proofs of formulae~\Ref{eii} are straightforward.

The proofs of formulae~\Ref{tT11} and \Ref{tT22} are similar. We will prove
formula \Ref{tT22}. Clearly, the weight subspace $\Wy[\>l\>]$ is spanned
by vectors of the form
\beq
\label{tT12}
\tT_{12}(u_1)\>\dots\>\tT_{12}(u_l)\cdot1\;=\;
y_0^l\,\prod_{i=1}^l\,\prod_{j=1}^{n-1}\,(u_i-y_j)
\eeq
with various $u_1,\dots,u_l$. So, it suffices to verify formula \Ref{tT22}
on such vectors.

Both the expression $\tT_{22}(u)\>\tT_{12}(u_1)\>\dots\>\tT_{12}(u_l)\cdot1$
and the right-hand side of formula \Ref{tT22} applied to
$\tT_{12}(u_1)\>\dots\>\tT_{12}(u_l)\cdot1$ are polynomials in $u$ of degree
$n$. Therefore, they are uniquely determined by their coefficients at $u^n$
and $u^{n-1}$, and the values at $n-1$ points $y_1,\dots,y_{n-1}$.

Proposition \ref{comm rel} and formulae~\Ref{Tab1}, \Ref{def of tilde-T},
\Ref{tT12} yield that
\be
\tT_{22}(u)\>\tT_{12}(u_1)\>\dots\>\tT_{12}(u_l)\cdot1\;=\;
\Bigl(u^n+\bigl(\>l-\sum_{i=1}^n z_i\bigr)\>u^{n-1}\Bigr)\,
y_0^l\,\prod_{i=1}^l\,\prod_{j=1}^{n-1}\,(u_i-y_j)\;+\;O(u^{n-2})
\ee
as $u\to\infty$, and
\vvn-.4>
\be
\bigl(\tT_{22}(u)\>\tT_{12}(u_1)\>\dots\>\tT_{12}(u_l)\cdot1\bigr)\big|_{u=y_j}
\;=\;d(y_j)\,y_0^l\,\prod_{i=1}^l\,\Bigl((u_i-y_j+1)\,
\prod_{j'\ne j}\,(u_i-y_j)\Bigr)\ ,
\vv-.6>
\ee
which proves the theorem.
\end{proof}

\subsubsection{}
\label{Skl cor}
{\bf Corollary.}\enspace
{\it We have}
\vvn-.3>
\begin{align*}
B(u, \bs{\tilde H})\ =\ \tT_{11}(u) + \tT_{22}(u) \ =
\ (2u + \sum_{i=1}^n (m_i-2z_i) + 2\sum_{j=1}^{n-1} y_j)\;
\prod_{j=1}^{n-1}\,(u-y_j) \ +{} &
\\[3pt]
{}\>+\ \sum_{j=1}^{n-1} \biggl(\,\prod_{j' \ne j}
\frac {u-y_{j'}}{y_j-y_{j'}}\>\biggr)
\left(a(y_j)\shift_{y_j}^{-1} + d(y_j)\shift_{y_j}\right)\, &\,\,.
\end{align*}

\subsection{Universal weight function}
\label{Universal weight function}
Let $\bs y = (y_0,\dots,y_{n-1})$. Recall that
$\bs a = (a_1,\alb\dots,a_l)$\>, \,$\bs h = (h_1,\dots,h_n)$ \,and
\,$p(x,\bs a) = x^l + a_1x^{l-1} + \dots + a_l$. Let
\be
\om(\bs y, \bs a) \ =\ y_0^l\,\prod_{j=1}^{n-1}\,p(y_j-1, \bs a)\ .
\ee
This element of $\Wy[\>l\>] \otimes \C[\bs a]\subset
\Wy[\>l\>] \otimes \C[\bs a, \bs h]$ is called {\it
the universal weight function.}

\medskip
A trivial but important property of the universal weight function is
given by the following lemma.

\subsubsection{}
\label{omega is nonzero}
{\bf Lemma.}
{\it Consider $\C^{l+n}$ with coordinates $\bs a, \bs h$. Then for
every $\T \in \C^{l+n}$, the vector $\om(\bs y, \bs a(\T))$ is
a nonzero vector of\, $\Wy[\>l\>]$}.
\qed

\medskip
Denote by $\om_D$ the projection of the universal weight function
$\om(\bs y, \bs a)$ to \,$\Wy[\>l\>]\otimes A_D\>=\,\W[\>l\>]\otimes A_D$.

\subsubsection{}
\label{thm on Bethe ansatz}
{\bf Theorem.}
{\it
Assume that the pair $\bigl((m_1,0),\dots, (m_n,0)\bigr)\>,\,l$ is separating.
Then for $s=1,\dots,n$, we have
\beq
\label{H equals h}
\tilde H_s\>\om_D\,=\,h_s\>\om_D\
\eeq
in \,$\Wy[\>l\>]\otimes A_D$. Moreover, we have}
\vvn.2>
\beq
\label{omega is sing}
\om_D\,\in\,\Sing\,\Wy[\>l\>]\otimes A_D\subset \Wy[\>l\>]\otimes A_D\ .
\vv.2>
\eeq

\subsubsection{}
\label{Homyap}
{\bf Corollary.}
{\it Let $\T$ be a point of the scheme \,$C_D$ considered as a set. Then
\vvn.3>
\beq
\label{omyapsing}
\om(\bs y, \bs a(\T))\,\in\,\Sing\,\Wy[\>l\>]\;.
\eeq
Moreover, for $s=1,\dots,n$, we have}
\beq
\label{H=h}
H_s\>\om(\bs y, \bs a(\T))\,=\,h_s(\T)\,\om(\bs y, \bs a(\T))\;.
\eeq

\begin{proof}[Proof of Corollary~\ref{Homyap}]
Let \,$\pi:\C[\bs a,\bs h]\to A_D$ \,be the canonical projection. A point
$\T\in C_D$ determines uniquely an algebra homomorphism $\Hat\T:A_D\to\C$,
such that $f(\T)=\Hat\T\bigl(\pi(f)\bigr)$ for any $f\in\C[\bs a,\bs h]$.
In particular,
\vvn-.3>
\beq
\label{omyap}
\om(\bs y, \bs a(\T))\,=\,(\id\otimes\Hat\T)(\om_D)\,.
\vv.3>
\eeq
Therefore, formulae~\Ref{omyapsing} and~\Ref{H=h} follow from
formulae~\Ref{omega is sing} and~\Ref{H equals h}, respectively.
\end{proof}

\subsubsection{}
{\bf Corollary.}
\label{linind}
{\it Let \,$\T_1,\dots\T_d$ be distinct points of the scheme \,$C_D$
considered as a set. Then the vectors
\,$\om(\bs y,\bs a(\T_1)),\dots,\om(\bs y,\bs a(\T_d))$ \,are linearly
independent.
}

\begin{proof}[Proof of Corollary~\ref{linind}]
The vector $\om(\bs y,\bs a(\T_j))$ is nonzero by Lemma~\ref{omega is nonzero}
and is an eigenvector of the operator $H_s$ with eigenvalue $h_s(\T_j)$
by formula~\Ref{H=h}. Moreover, the collections of eigenvalues
$\bs h(\T_1),\dots,\alb\bs h(\T_d)$ are distinct, because a point $\T\in C_D$
is uniquely determined by its coordinates $\bs h(\T)$ by
Theorem~\ref{thm-conjecture 2}. The corollary is proved.
\end{proof}

\begin{proof}
[Proof of Theorem~\ref{thm on Bethe ansatz}]
To prove formula
\Ref{H equals h} it is enough to show that the polynomial
$\bigl(B(u, \bs{\tilde H}) - B(u, \bs h)\bigr)\,\om(\bs y, \bs a)$
projects to zero in $\C[u]\otimes \Wy[\>l\>]\otimes A_D$.
Let
\bea
\Bc(u,y_1,\dots,y_{n-1}, \bs h)\ =
\ \sum_{j=1}^{n-1}\, B(y_j, \bs h)
\prod_{j' \ne j} \frac {u-y_{j'}}{y_j-y_{j'}}\ .
\eea
For $j=1,\dots,n$, we have
\,$\Bc(y_j,y_1,\dots,y_{n-1}, \bs h)\,=\,B(y_j,\bs h)$ \,and
$\Bc(u,y_1,\dots,y_{n-1}, \bs h)$ is a polynomial in $u$ of degree
$\leq n-2$. Hence
\vvn-.2>
\be
B(u,\bs h) - \Bc(u,y_1,\dots,y_{n-1}, \bs h)
\ = \ \Bigl(2u+h_1+ 2\sum_{j=1}^{n-1} y_j\Bigr)\,\prod_{j=1}^{n-1}(u-y_j)\ .
\vv-.4>
\ee
We have
\begin{align*}
\bigl(B & (u,\bs{\tilde H}) - B(u,\bs h) +
\Bc(u,y_1,\dots,y_{n-1}, \bs h) -
\Bc(u,y_1,\dots,y_{n-1}, \bs h))\om(\bs y, \bs a\bigr)\ ={}
\\[4pt]
& \biggl(\!\Bigl(-\>h_1 + \sum_{i=1}^n (m_i-2z_i)\Bigr)\,
\prod_{j=1}^{n-1}\,(u-y_j)\biggr)\>\om(\bs y, \bs a) \ +{}
\\[4pt]
& \hphantom{\biggl(}\sum_{j=1}^{n-1}\,y_0^l\,\biggl(
\,\prod_{j' \ne j}\>\frac {u-y_{j'}}{y_j-y_{j'}}\;p(y_{j'}-1,\bs a)\biggr)
\!\left(a(y_j)\shift_{y_j}^{-2} - B(y_j, \bs h)
\shift_{y_j}^{-1}+ d(y_j)\right)p(y_j,\bs a)\ .
\end{align*}
Clearly all terms in the right-hand side of this formula project to zero
in $\C[u]\otimes \Wy[\>l\>]\otimes A_D$.
Hence, formula \Ref{H equals h} is proved.

\medskip
The proof of formula \Ref{omega is sing} is based on the following
lemma.

\medskip
\noindent
{\bf Lemma.}
{\it We have $e_{21}e_{12}\>\om_D\,=\,0$. }

\begin{proof}
{}From the formula for the quantum determinant we have
\vvn.3>
\beq
\label{hor}
\tT_{12}(u)\tT_{21}(u-1)\om(\bs y,\bs a)\ =
\ \bigl(\tT_{11}(u)\tT_{22}(u-1)\ -\ a(u)d(u-1)\bigr)\,\om(\bs y,\bs a)\ ,
\vv.3>
\eeq
where $\tT_{12}(u)\tT_{21}(u-1)\,=\,e_{21}e_{12}u^{2n-2} + O(u^{2n-3})$.
Therefore, our goal is to calculate the coefficient of $u^{2n-2}$ in
the right-hand side. We have
\be
T_{11}(u)T_{22}(u-1)\ =\ 1\,+\,\frac{e_{11}+e_{22}}{u}\,+\,
\frac{e_{22}(e_{11}+1) + T^{\{2\}}_{11}+T^{\{2\}}_{22}}{u^2}\,+\,O(u^{-3})\ .
\ee
Hence
\be
T_{11}(u)T_{22}(u-1) - T_{11}(u) - T_{22}(u) + 1\ =
\ \frac{e_{22}(e_{11}+1)}{u^2}\,+\,O(u^{-3})\
\vv-.3>
\ee
and
\vvn.2>
\be
\tT_{11}(u)\tT_{22}(u-1) - B(u, \bs{\tilde H})\>d(u-1) + d(u)d(u-1)\ =
\ e_{22}(e_{11}+1) u^{2n-2} + O(u^{2n-3})\ .
\vv.2>
\ee
Thus the right-hand side of \Ref{hor} equals
\vvn.2>
\be
\bigl(B(u, \bs{\tilde H}) - a(u) - d(u)\bigr)d(u-1)\,\om(\bs y,\bs a)
+ e_{22}(e_{11}+1) u^{2n-2}\om(\bs y,\bs a) + O(u^{2n-3})\ .
\vv.2>
\ee
Here
\vvn->
\begin{align*}
e_{22}(e_{11}+1)\,\om(\bs y,\bs a) \ &{}=
\ l \>\Bigl(\sum_{i=1}^n m_i - l+ 1\Bigr) \,\om(\bs y,\bs a)\ ,
\\
B(u, \bs{\tilde H})\, \om(\bs y,\bs a) \ &{}=
\ B(u, \bs{h}) \,\om(\bs y,\bs a)\ ,
\end{align*}
\be
a(u) + d(u) = 2u^n - \sum_{s=1}^n\,(2z_s - m_s) u^{n-1}
+ \!\! \sum_{1\leq i< j\leq n}
\bigl(z_iz_j + (z_i-m_i)(z_j-m_j)\bigr)\,u^{n-2} + {}\dots{}
\vv.2>
\ee
and \;$B(u, \bs{h})\,=\,2u^n + h_1 u^{n-1} + h_2 u^{n-2} +{}\dots{}$\>.
Therefore, the right-hand side of \Ref{hor} equals
\vvn-.2>
\begin{align*}
& \Bigl(h_1 + \sum_{s=1}^n (2z_s - m_s)\Bigr)\,u^{n-1} d(u-1)\,
\om(\bs y,\bs a)\,+{}
\\
& \Bigl(h_2 +l \Bigl(\sum_{i=1}^n m_i - l+ 1\Bigr)
- \!\!\sum_{1\leq i< j\leq n}\!\bigl(z_iz_j + (z_i-m_i)(z_j-m_j)\bigr)\!\Bigr)
\>u^{2n-2}\,\om(\bs y,\bs a) + O(u^{2n-3})\ .
\end{align*}
Clearly the first two terms of this expression project to zero in
$\C[u]\otimes \Wy[\>l\>]\otimes A_D$. This proves the lemma.
\end{proof}

In order to deduce formula \Ref{omega is sing} from the lemma,
it is enough to
notice that the operator $e_{21}$ is injective,\ in variables
$\bs y$ it is the operator of multiplication by $y_0$.
Therefore, $e_{12} \,\om_D\, =\, 0$.
Theorem \ref{thm on Bethe ansatz} is proved.
\end{proof}

\section{Multiplication in algebra $A_D$ and Bethe algebra $A_W$}
\label{Multiplication Mult}

\subsection{Multiplication in $A_D$}

By Theorem \ref{thm const wrt z},
the scheme $C_D$ considered as a set is finite,
and the algebra $A_D$ is the direct sum of local algebras,
\be
A_D\ =\ \oplus_{\T}\ A_{\T, D}\
\ee
corresponding to points $\T$ of the set $C_D$.

The local algebra $\AD$ may be defined as the quotient of the algebra of germs
at $\T$ of holomorphic functions in $\bs a, \bs h$ modulo the ideal $I_{\T,D}$
generated by all functions $q_{1},q_2,\dots,q_{l+n}$.

The local algebra $\AD$ contains the maximal ideal $\mg_\T$
generated by germs which are zero at $\T$.

For $f\in A_D$, denote by $L_f$ the linear operator $A_D \to A_D$\>,
\,$g\mapsto fg,$\ of multiplication by $f$. Consider the dual space
\be
A_D^*\ =\ \oplus_\T\, \AD^*
\vv.2>
\ee
and the dual operators $L_f^* : A_D^* \to A_D^*$.

Every summand $\AD^*$ contains the distinguished one-dimensional subspace
$\mg_\T^\perp$ which is the annihilator of $\mg_\T$.

\subsubsection{}
\label{lemma on dual operators}
{\bf Lemma} \cite{MTV3}.
{\it
\begin{enumerate}
\item[(i)]
For any point $\T$ of the scheme $C_D$ considered as a set and any $f\in A_D$,
we have $L_f^* (\mg_\T^\perp) \subset \mg_\T^\perp$.

\item[(ii)] For any point $\T$ of the scheme $C_D$ considered as a set,
if \;$W\subset \AD^*$ is a nonzero vector subspace invariant with respect to
all operators $L_f^*$, $f\in A_D$, then \,$W$ contains $\mg_\T^\perp$.
\end{enumerate}
}

\begin{proof} For any $f \in \mg_\T$ we have
$L_f^*(\mg_\T^\perp) = 0$. This proves part (i).

To prove part (ii) we
consider the filtration of $\AD$ by powers of the
maximal ideal,
\be
\AD \supset \mg_\T \supset \mg_\T^2\supset\dots\supset \{0\}\ .
\ee
We consider a linear basis $\{f_{a,b}\}$ of $\AD$, \,$a=0,1,\dots{}$,
\,$b=1,2,\dots{}$, which agrees with this filtration. Namely, we assume that
for every $i$, the subset of all vectors $f_{a,b}$ with $a\geq i$ is a basis
of $\mg^i_\T$\,.

Since dim $\AD/\mg_\T = 1$, there is only one basis vector with $a=0$ and
we also assume that this vector $f_{0,1}$ is the image of $1$ in $\AD$\,.

Let $\{ f^{a,b}\}$ denote the dual basis of $\AD^*$. Then
the vector $f^{0,1}$ generates $\mg_\T^\perp$.

Let $w = \sum_{a,b} c_{a,b} f^{a,b}$ be a nonzero vector in $W$. Let
$a_0$ be the maximum value of $a$ such that there exists $b$ with a nonzero
$c_{a,b}$. Let $b_0$ be such that $c_{a_0,b_0}$ is nonzero.
Then it is easy to see that $L^*_{f_{a_0,b_0}} w\,=\,c_{a_0,b_0} f^{0,1}$.
Hence $W$ contains $\mg_\T^\perp$.
\end{proof}

\subsection{Linear map $\tau : A_D^* \to \Sing\,\W[\>l\>]$}
Let $f_1,\dots, f_\mu$ be a basis of $A_D$
considered as a vector space over $\C$. Write
\vvn.2>
\beq
\label{basis in M}
\om_D\ =\ \sum_i v_i\otimes f_i\qquad
{\rm with} \qquad v_i \in \Sing\,\Wy[\>l\>] = \snash \ .
\eeq
Denote by $V \subset\snash$ the vector subspace spanned by
$v_1,\dots,v_\mu$. Define the linear map
\vvn-.3>
\beq
\label{map psi}
\tau\ :\ A_D^*\ \to\ \snash\ , \qquad
g\ \mapsto\ g(\om_D) = \sum_i\ g(f_i)\,v_i\ .
\eeq
Clearly, $V$ is the image of $\tau$.

\subsubsection{}
\label{lemma april 1}
{\bf Lemma.}
{\it Let $\T$ be a point of $C_D$ considered as a set. Let
\vvn.3>
\be
\om(\bs y, \bs a(\T)) \in \Wy[\>l\>] = \W[\>l\>]
\vv.3>
\ee
be the value of the universal weight function at $\T$.
Then the vector $\om(\bs y, \bs a(\T))$ belongs to the image of $\tau$.}

\begin{proof}
The statement follows from formula~\Ref{omyap}.
\end{proof}

Let $\psi_{DW}:A_D \to A_W$ be the epimorphism defined in
Theorem~\ref{thm D to B}.

\subsubsection{}
\label{lemma 1}
{\bf Lemma}
{\it Assume that the pair $\bigl((m_1,0),\dots, (m_n,0)\bigr)\>,\,l$ is
separating. Then for any $f\in A_D$ and $g\in A_D^*$, we have
\,$\tau (L^*_f(g)) = \psi_{DW}(f) (\tau(g))$\>.}

\medskip
In other words, the map $\tau$ intertwines the action of the algebra
of multiplication operators $L^*_f$ on $A_D^*$ and the action on
the Bethe algebra on $\snash$.

\begin{proof}
The algebra $A_D$ is generated by $h_1,\dots,h_n$. It is enough to prove
that for any $s$ we have $\tau (L^*_{h_s} (g))\, =\, H_s (\tau(g))$.
But $\tau(L^*_{h_s} (g)) = \sum_i g(h_s f_i) v_i =
g\bigl(\sum_i \,v_i \otimes h_sf_i\bigr) =
g\bigl(\sum_i \,H_s v_i \otimes f_i \bigr) = H_s (\tau(g))$.
\end{proof}

\subsubsection{}
\label{corollary 1}
{\bf Corollary.}
{\it
The vector subspace $V\subset\snash$
is invariant with respect to the action of the Bethe algebra $A_W$
and the kernel of $\tau$ is a subspace of $A_D^*$, invariant with respect to
multiplication operators $L^*_f,\,f\in A_D$.
}

\subsection{First main theorem}
\subsubsection{}
\label{1st main thm}
{\bf Theorem.}
{\it
Assume that the pair $\bigl((m_1,0),\dots, (m_n,0)\bigr)\>,\,l$ is separating.
Then the image of $\tau$ is $\Sing\,\W[\>l\>]$ and the kernel of $\tau$ is
zero.}

\subsubsection{}
\label{1st main cor}
{\bf Corollary.}
{\it
The map $\tau$ identifies the action of operators $L_f^*$, $f\in A_D$,
on $A_D^*$ and the action of the Bethe algebra on $\Sing\,\W[\>l\>]$.
Hence the epimorphism $\psi_{DW}\, : \, A_D\, \to\, A_W$\ is an isomorphism.}

\goodbreak
\begin{proof}[Proof of Theorem \ref{1st main thm}]
First we will show that $\tau$ is an epimorphism for generic $\bs z$.

Let $d_l\,=\,\dim\,\snash$. Corollary \ref{cor C_D has many points} says
that for generic $\bs z$ there exists $d_l$ distinct points
$\T_1,\dots,\T_{d_l}$ in $C_D$. By Corollary~\Ref{linind}, the vectors
$\om(\bs y, \bs a(\T_1))$, \dots , $\om(\bs y, \bs a(\T_{d_l}))$
are linearly independent and hence form a basis for $\snash$. Therefore,
$\tau$ is an epimorphism for generic $\bs z$ by Lemma \ref{lemma april 1}.

By Theorem \ref{thm const wrt z} and Lemma \ref{singular}, dimensions of $A_D$
and $\snash$ do not depend on $\bs z$. Hence $\dim\,A_D \geq \dim\,\snash$ for
all $z_1,\dots, z_n$. Therefore, to prove Theorem \ref{1st main thm}
it remains to prove that $\tau$ has zero kernel.

Denote the kernel of $\tau$ by $K$. Let $A_D = \oplus_\T \AD$ be
the decomposition into the direct sum of local algebras.
Since $K$ is invariant with respect to multiplication operators, we have that
$K \,= \,\oplus_\T\, K\cap \AD^*$, and for every $\T$\>, the vector subspace
$K\cap \AD^*$ is invariant with respect to multiplication operators.
By Lemma \ref{lemma on dual operators}, if $K\cap \AD^*$ is nonzero,
then $K\cap \AD^*$ contains the one-dimensional subspace $\mg_\T^\perp$.

Let $\{f_{a,b}\}$ be the basis of $\AD$ constructed in the proof of Lemma
\ref{lemma on dual operators}, and let $\{f^{a,b}\}$ be the dual basis of
$\AD^*$. Then the vector $f^{0,1}$ generates $\mg_\T^\perp$. By definition
of $\tau$, the vector $\tau(f^{0,1})$ is equal to the value of the universal
weight function at $\T$. By Lemma \ref{omega is nonzero}, this value is
nonzero and that contradicts the assumption that $f^{0,1}$ lies in
the kernel of $\tau$.
\end{proof}

\subsection{Grothendieck bilinear form on $A_D$}
\label{Grothendieck bilinear form on A-D}

Realize the algebra $A_D$ as $\C[\bs h]/I'_D$,
\,where $I'_D$ is the ideal generated by $n$ polynomials
$q_{1}(\bs h)\>,\,q_2(\bs h)\>,\,q_j(\bs a(\bs h),\bs h)$\>,
\,$j= l+3,\dots,\alb l+n$, see~\Ref{sec 2nd description}.

Let $\varrho: A_D \to \C$, be the Grothendieck residue,
\vvn.2>
\be
f \ \mapsto \ \frac 1{(2\pi i)^n}\,\Res_{C_D}
\frac{ f}{q_{1}(\bs h) q_2(\bs h)\prod_{j=l+3}^{l+n}\,
q_j(\bs a(\bs h),\bs h)}\ .
\vv.3>
\ee
Let $(\,,\,)_D$ be the Grothendieck symmetric bilinear form on $A_D$
defined by the rule
\vvn.2>
\be
( f,\, g)_D\ =\ \varrho (f g)\ .
\vv.2>
\ee
The Grothendieck bilinear form is nondegenerate.

The form $(\,,\,)_D$ determines a linear isomorphism $\phi:A_D \to A_D^*$,
\;$f \mapsto (f,\,\cdot)_D$.

\subsubsection{}
\label{lem on intertw}
{\bf Lemma.}
{\it
The isomorphism $\phi$ intertwines the operators $L_f$ and $L^*_f$
for any $f\in A_D$. }

\begin{proof} For $g\in A_D$ we have $\phi (L_f(g))=
\phi (fg) = (fg,\cdot)_D = (g, f\cdot)_D = L_f^*((g,\cdot)_D) =
L_f^*\phi(g)$.
\end{proof}

\subsubsection{}
\label{intersting cor}
{\bf Corollary.} {\it
Assume that the pair $\bigl((m_1,0),\dots, (m_n,0)\bigr)\>,\,l$ is separating.
Then the composition $\tau\phi\,:\, A_D \to \snash$ is a linear isomorphism
which intertwines the algebra of multiplication operators on $A_D$ and
the action of the Bethe algebra $A_W$ on $\snash$.
}

\section{Algebra $A_G$} % Three more algebras
\label{Three more algebras}

\subsection{New conditions on $(m_1,0)\ ,\ \dots\ ,\ (m_n,0),\,l$}
\label{sec new conditiong}
In the remainder of the paper we assume that
\vvn-.3>
\be
\bs \La\ =\ (\La^{(1)}, \dots, \La^{(n)}) \ =\
\bigl((m_1,0)\ ,\ \dots\ ,\ (m_n,0)\bigr)
\vv.2>
\ee
is a collection of dominant integral
$\glt$-weights, that is, $m_s \in \Z_{\geq 0}$ for
$s = 1,\dots, n$.

We assume that $l \in \Z_{\geq 0}$ is such that the weight
$\bigl(\>\sum_{s=1}^nm_s - l,\,l\>\bigr)$ is dominant integral,
that is, $\sum_{s=1}^nm_s - l \geq l$. This assumption implies that
the pair $\bigl((m_1,0)\ ,\ \dots\ ,\ (m_n,0)\bigr)\>,\,l$ is separating.

Let \,$\tilde l\,=\>\sum_{s=1}^n m_s\,+1 - l$. We have \,$\tilde l > l$.

\subsection{Wronskian}
\label{sec Wronskian}
The (discrete) Wronskian of polynomials $f,g \in \C[u]$ is
\vvn.2>
the polynomial
\be
\Wr\,(f(u),g(u))\ =\ f(u)g(u-1) - f(u-1)g(u)\ .
\ee

\subsubsection{}
\label{lemma on Wr}
{\bf Lemma.}
{\it
Let $f,g, \Bc\in \C[u]$. Assume that $f, g$ are
monic polynomials of degrees $l, \tilde l$, respectively,
that lie in the kernel of the difference operator
\vvn.4>
\be
d(u)\ -\ \Bc(u)\, \shift^{-1}\ +\ a(u) \, \shift^{-2}\ .
\ee
Then}
\vvn-.4>
\be
\Wr\,(f(u),g(u))\ =\ (l - \tilde l)
\prod_{s=1}^n \prod_{j=1}^{m_s}\,(u-z_s + j)\ .
\vv-.5>
\ee
\begin{proof}
Let $\mathcal C(u)=\Wr\,(f(u),g(u))$. Then the top coefficient of
$\mathcal C(u)$ equals $l-\tilde l$, and
\vvn.1>
\be
\frac{\mathcal C(u)}{\mathcal C(u-1)}\;=\;\frac{a(u)}{d(u)}\ ,
\vv.1>
\ee
which determines the polynomial $\mathcal C(u)$ uniquely.
\end{proof}

\subsubsection{}
\label{lemma 2 on Wr}
{\bf Lemma.}
{\it Let $f,g \in \C[u]$, $z\in \C$, $m \in \Z_{> 0}$.
Assume that
$f(z-j)=0$ for $j=1,\dots,m+1$.
Then the polynomial $\Wr\,(f(u),g(u))$ is equal to zero at
$u = z-j$, $j=1, \dots , m$, and the polynomial $f(u) g(u-2) - f(u-2) g(u)$
is equal to zero at $u = z-j$, $j=1, \dots , m-1$.
}
\qed

\subsubsection{}
\label{lemma 3 on Wr}
{\bf Lemma.}
{\it Let $f,g, \mathcal C \in \C[u]$, $z\in \C$,
and $\Wr\,(f(u),g(u))\,=\,\mathcal C(u)$.
\begin{enumerate}
\item[(i)] If $\mathcal C(z)\ne 0$ and $f(z-1)=0$, then $g(z-1)\ne 0$.
\item[(ii)] If $\mathcal C(z)\ne 0$ and $f(z)=0$, then $f(z-1)\ne 0$.
\qed
\end{enumerate}}

\subsubsection{}
\label{lemma 5 on Wr}
{\bf Lemma.}
{\it Let $f,g \in \C[u]$, $z\in \C$. Then}
\vvn.5>\\
\rightline{$\displaystyle
\Wr\,((u-z)f(u),(u-z)g(u))\ = \ (u-z)(u-z-1)\,\Wr\,(f(u),g(u))\ .$
\hfil\llap{$\square$}}

\subsection{Intersection of Schubert cycles $C_G$}
\label{sec Schubert cells and cycles}

Let $d$ be a sufficiently large natural number
with respect to the numbers $m_1,\dots, m_n$ considered in Section
\ref{sec new conditiong}. Let
$\C_d[u]$ be the vector subspace in $\C[u]$ of polynomials
of degree not greater than $d$.

\subsubsection{}
\label{sec def of schub}
Denote by $G$ the Grassmannian of all two-dimensional vector subspaces in
$\C_d[u]$.

Let $\F = \{ 0 = F_{d+1}
\subset F_d \subset \dots \subset F_1 \subset
F_0 = \C_d[u]\}$ be a complete flag and $\La = (a,b)$ a $\glt$
dominant integral weight such that $d \geq a\geq b \geq 0$ and $a,b \in \Z$.
Define a Schubert cell $C^o_{\F, \La} \subset G$ to be
the set of all two-dimensional subspaces
$V \subset \C_d[u]$ having a basis $f, g$ such that
\vvn-.5>
\be
f \in F_{a+1} - F_{a+2} \qquad {\rm and} \qquad g \in F_b - F_{b+1}\ .
\ee
Define a Schubert cycle \,$C_{\F, \La} \subset G$ as the closure of
the Schubert cell \,$C^o_{\F, \La}$.

For $z \in \Z$ and $i \in \Z_{>0}$, set
\vvn-.6>
\be
\varphi_i(u,z)\ =\ \prod_{j=1}^i(u-z+j)\ .
\vv-.4>
\ee
Introduce a complete flag in $\C_d[u]$\,:
\be
\F(z)\ = \
\{ 0 = F_{d+1}(z) \subset F_d(z) \subset \dots \subset F_1(z) \subset
F_0(z) = \C_d[u]\}\ ,
\ee
where $F_i(z)$ consists of all polynomials divisible by $\varphi_i(u,z)$.

Introduce the complete flag in $\C_d[u]$ associated with infinity:
\be
\F(\infty)\ = \
\{ 0 = F_{d+1}(\infty) \subset F_d(\infty)
\subset \dots \subset F_1(\infty) \subset
F_0(\infty) = \C_d[u]\}\ ,
\ee
where $F_i(\infty)$ consists of all polynomials of degree $\leq d-i$.

\smallskip
We consider the Schubert cells \,$C^o_{\F(z_s), \La^{(s)}} \subset G$,
\,$s = 1, \dots, n$, \,where $\La^{(s)}=(m_s,0)$, and the Schubert cell
\,$C^o_{\F(\infty), \La^{(\infty)}} \subset G$,
where $\La^{(\infty)} = (d - l, d- \tilde l - 1)$.
The cell \,$C^o_{\F(z_s), \La^{(s)}}$ is the
set of all two-dimensional subspaces $V \subset \C_d[u]$ having a basis $f, g$
such that
\be
g(z_s-1) \ne 0\ ,\qquad f(z_s-m_s-2) \ne 0\ ,\qquad f(z_s-j)=0
\quad {\rm for}\ j=1,\dots,m_s+1 \ ,
\ee
and the cell $C^o_{\F(\infty), \La^{(\infty)}}$ is the set of all
two-dimensional subspaces $V \subset \C_d[x]$ having a basis $f, g$
\,such that \,$\deg\, f = l$ \,and \,$\deg\, g = \tilde l$.

Consider the (scheme-theoretic) intersection
\vvn.2>
\beq
\label{CG}
C_G\ =\ C_{\F(\infty), \La^{(\infty)}} \,\tbigcap\,\bigl(\,\cap_{s=1}^n\,
C_{\F(z_s), \La^{(s)}}\,\bigr)\
\vv.2>
\eeq
of the corresponding Schubert cycles.
Denote by $A_G$ the algebra of functions on $C_G$.

\subsubsection{}
\label{lemma on open intersection}
{\bf Lemma.}
{\it Let \,$z_i-z_j \notin \Z$ \,for \,$i\ne j$. Then }
\vvn.2>
\be
C_G\,=\,C^o_{\F(\infty), \La^{(\infty)}}\,\tbigcap\,\bigl(\,\cap_{s=1}^n\,
C^0_{\F(z_s), \La^{(s)}}\,\bigr)\
\ee
as sets.

\begin{proof}
Let \,$V$ be a point of
$ C_{\F(\infty), \La^{(\infty)}} \, \tbigcap\,\bigl( \cap_{s=1}^n\,
C_{\F(z_s), \La^{(s)}}\bigr)$\>. Let $f, g$ be a monic basis of $V$,
\vvn.1>
such that $\deg f \leq l$ \,and \,$\deg\>g \leq \tilde l$.
\vvn.2>
Then \,$\deg \Wr\,(f(u),g(u)) \leq l+\tilde l - 1$.
On the other hand, the polynomial $\Wr\,(f(u),g(u))$ is divisible by
\vvn.1>
$\prod_{s=1}^n \prod_{j=1}^{m_s}\,(u-z_s + j)$ by Lemma~\ref{lemma 2 on Wr}.
Since \,$\sum_{s=1}^n m_s = \tilde l + l -1$, we conclude that
\;$\deg f = l$, \;$\deg\>g = \tilde l$,
\;$V$ is a point of $ C^o_{\F(\infty), \La^{(\infty)}}$\>, \;and
\vvn-.6>
\be
\Wr\,(f(u),g(u))\,=\,(l-\tilde l)\,\prod_{s=1}^n\prod_{j=1}^{m_s}\,(u-z_s+j)\,.
\ee
Since a suitable linear
combination of $f$ and $g$ is divisible by \,$\varphi_{m_s+1}(u,z_s)$,
the subspace $V$ is a point of
$C^o_{\F(z_s), \La^{(s)}}$ by Lemma \ref{lemma 3 on Wr}.
\end{proof}

\subsubsection{}
Let $V$ be a point of $C_G$ considered as a set.
Then there exists a unique basis $f, g$ of $V$ such that
\begin{align*}
f(u)\ &{}=\ u^l + f_1u^{l-1}+ \dots + f_l\ ,
\\[3pt]
g(u)\ &{}=\ u^{\tilde l} + {g}_1u^{\tilde l-1}+ \dots +
g_{\tilde l - l -1} u^{l+1} +
g_{\tilde l - l +1}u^{l-1}+ \dots +g_{\tilde l}\
\end{align*}
for suitable complex numbers $f_1,\dots,f_l,\,
g_1, \dots, g_{\tilde l - l -1}, g_{\tilde l - l +1}, \dots, g_{\tilde l}$.

\subsubsection{}
\label{lemma on diff eqn}
{\bf Lemma.}
{\it
Let $z_i-z_j \notin \Z$ for \,$i\ne j$. Then all polynomials of
the subspace $V$ are annihilated by the difference operator
\vvn.2>
\be
\D_V\,=\,d(u)\,-\,\Bc_V(u)\, \shift^{-1}\>+\,a(u) \, \shift^{-2}\ ,
\vv-.5>
\ee
where
\vvn-.5>
\be
\Bc_V(u)\,=\,\frac 1{\tilde l-l}\,\bigl(g(u)f(u-2) - g(u-2)f(u)\bigr)\,
\prod_{s=1}^n \prod_{j=1}^{m_s-1}\,(u-z_s + j)^{-1}
\ee
is a polynomial of degree $n$.
}

\begin{proof}
Let $W(u)=\Wr\,\bigl(f(u),g(u)\bigr)$.
It is straightforward to see that all polynomials of
the subspace $V$ are annihilated by the difference operator
\vvn.2>
\beq
\label{WBW}
W(u-1)\,-\,\bigl(g(u)f(u-2) - g(u-2)f(u)\bigr)\,\shift^{-1}\>+\,
W(u)\,\shift^{-2}\,.
\vv-.1>
\eeq
Since
\vvn-.5>
\be
W(u)\,=\,(l-\tilde l)\,\prod_{s=1}^n\prod_{j=1}^{m_s}\,(u-z_s+j)\,,
\ee
see the proof of Lemma~\ref{lemma on open intersection}, and
all coefficients of the difference operator~\Ref{WBW} are divisible by
$\prod_{s=1}^n \prod_{j=1}^{m_s-1}\,(u-z_s+j)$ by Lemma~\ref{lemma 2 on Wr},
the statement follows.
\end{proof}

Write
\vvn-.6>
\be
\Bc_V(u)\,=\,2u^n + h_1u^{n-1} + \dots + h_n\,.
\ee
Recall that the scheme $C_D$ is defined in Section~\ref{subsub h's}.

\subsubsection{}
\label{CGCD}
{\bf Corollary of Lemma \ref{lemma on diff eqn}.}
{\it Consider the schemes $C_D$ and $C_G$ as sets. Then the
assignment \,$V\mapsto(f_1,\dots,f_l, h_1,\dots,h_n)\in\C^{l+n}$
defines an injective map of sets\\
\,$C_G\to C_D$. }

\subsubsection{}
\label{thm inverse on Schub}
{\bf Theorem.}
{\it
Let \,$z_i-z_j \notin \Z$ \,for \,$i\ne j$.
Assume that $V \in G$ has a basis $f,g$ such that
$\deg f = l$ and $\deg g = \tilde l$, and $V$ is annihilated by
a difference operator of the form
\be
d(u)\,-\,\Bc(u)\, \shift^{-1}\>+\,a(u) \, \shift^{-2}\ ,
\vv.1>
\ee
where $\Bc(u)$ is a polynomial. Then $V$ is a point of \,$C_G$\,.}

\medskip\noindent
The proof is similar to the proof of Theorem~7.2 in \cite{MTV2}.
\qed

\subsection{Algebra $A_G$}
\label{sectAD}
\subsubsection{}
\vsk-.4>
\label{lemma on 0 dim}
{\bf Lemma.}
{\it
Let \,$z_i-z_j \notin \Z$ \,for \,$i\ne j$. Then
$A_G$ considered as a vector space is finite-dimensional.
Moreover, this dimension does not depend on $\bs z$.
}

\begin{proof}
The claim follows from Corollary~\ref{CGCD} and the reasoning similar to
the proof of Theorem~\ref{thm const wrt z}.
\end{proof}

Under conditions of Lemma \ref{lemma on 0 dim},
the dimension of $A_G$ as a vector space is given by Schubert calculus.
Namely, let $\bs \La = (\La^{(1)}, \dots , \La^{(n)})$ be the
collection of $\glt$-highest weights, where $\La^{(s)}=(m_s,0)$. Denote by
\vvn-.2>
\be
L_{\bs \La}\,=\,L_{\La^{(1)}}\otimes \dots \otimes L_{\La^{(n)}}
\vv.3>
\ee
the tensor product of irreducible $\glt$-modules with highest weights
$\La^{(1)}, \dots , \La^{(n)}$, respectively. Let
$\Sing \,L_{\bs\La}[\>l\>]$ be
the subspace of $L_{\bs\La}$ of
$\glt$-singular vectors of weight $(\sum_{s=1}^n m_s - l, l)$.
Then by Schubert calculus,
\beq
\label{dimAG}
\dim\,A_G\ =\ \dim\,\Sing\,L_{\bs \La}[\>l\>]\;,
\vv-.2>
\eeq
see \cite{Fu}.

\subsection{Presentation of algebra $A_G$}
\label{sec Presentation of algebra A_G}

If $z_i-z_j \notin \Z$ for $i\ne j$, we
shall use the following presentation of the algebra $A_G$.

Let
\vvn-.5>
\be
\tilde{\bs a} = (\tilde a_1, \dots, \tilde a_{\tilde l - l -1},
\tilde a_{\tilde l - l +1}, \dots, \tilde a_{\tilde l})\ .
\ee
Consider the space $\C^{\tilde l + l + n - 1}$ with coordinates
$\tilde{\bs a}, \bs a, \bs h$, cf. Section \ref{subsub h's}.

Denote by $\qat$ the following polynomial in $u$ depending on parameters
$\tilde{\bs a}$,
\be
\qat \ =
\ u^{\tilde l} + \tilde {a}_1u^{\tilde l-1}+ \dots +
\tilde a_{\tilde l - l -1} u^{l+1} +
\tilde a_{\tilde l - l +1}u^{l-1}+ \dots +\tilde a_{\tilde l}\ .
\ee
Recall that $\qa = u^l + a_1u^{l-1}+ \dots + a_l$ and
$B(u, \bs h) = 2u^n + h_1u^{n-1} + \dots + h_n$.

Let us write
\begin{align*}
& \Wr\,(\qat, \qa)\ = \ (\tilde l - l) u^{\tilde l + l - 1} +
w_1(\tilde {\bs a}, \bs a)u^{\tilde l + l - 2} + \dots +
w_{\tilde l + l -1}(\tilde {\bs a}, \bs a)\ ,
\\[4pt]
& \tilde p(u, \tilde{\bs a}) p(u-2, {\bs a}) -
\tilde p(u-2, \tilde{\bs a})p(u, {\bs a}) \ ={}
\\[2pt]
& \hphantom{\Wr\,(\qat, \qa)\ ={}\!}
2 (\tilde l - l) u^{\tilde l + l - 1} +
\hat w_1(\tilde {\bs a}, \bs a)u^{\tilde l + l - 2} + \dots +
\hat w_{\tilde l + l -1}(\tilde {\bs a}, \bs a)
\end{align*}
for suitable polynomials $w_1,\dots,w_{\tilde l + l -1}$,
$\hat w_1,\dots, \hat w_{\tilde l + l -1}$ in variables
$\tilde{\bs a}, \bs a$. Let us write
\begin{align*}
(\tilde l - l)\,\prod_{s=1}^n\,\prod_{j=1}^{m_s}\,(u-z_s+ j)\ ={}&
\ (\tilde l - l) u^{\tilde l + l - 1} +
c_1u^{\tilde l + l - 2} + \dots + c_{\tilde l + l -1}\ ,
\\
(\tilde l - l)\,B(u,\bs h) \,
\prod_{s=1}^n\prod_{j=1}^{m_s-1} (u-z_s+ j)\ ={}&
\ 2 (\tilde l - l) u^{\tilde l + l - 1} +
\hat c_1(\bs h)u^{\tilde l + l - 2} + \dots +
\hat c_{\tilde l + l -1}(\bs h)\ ,
\end{align*}
for suitable numbers $c_1,\dots,c_{\tilde l + l -1}$ and polynomials
$\hat c_1,\dots, \hat c_{\tilde l + l -1}$ in variables $ \bs h$.

Denote by $I_G$ the ideal in $\C[\tilde{\bs a}, \bs a, \bs h]$ generated by
$2(\tilde l + l -1)$ polynomials
\beq
\label{eqns for C_G}
w_i(\tilde{\bs a}, \bs a) - c_i\ , \qquad
\hat w_i(\tilde{\bs a}, \bs a) - \hat c_i(\bs h)\ ,
\qquad i = 1, \dots , \tilde l + l -1\ .
\eeq

\subsubsection{}
\label{lemma on new presentation}
{\bf Lemma.}
{\it
Let \,$z_i-z_j \notin \Z$ \,for \,$i\ne j$. Then}
\be
A_G\ = \ \C[\tilde{\bs a}, \bs a, \bs h] / I_G\ .
\ee

\begin{proof}
The scheme defined by the ideal $I_G$ consists of points $\T$ such that
\begin{align*}
\Wr\,(\tilde p(u, \tilde{\bs a}(\T)), p(u, {\bs a}(\T)))\ &{}=
\ (\tilde l - l)\,\prod_{s=1}^n\prod_{j=1}^{m_s} (u-z_s+ j)\ ,
\\
\tilde p(u, \tilde{\bs a}) p(u-2, {\bs a}) -
\tilde p(u-2, \tilde{\bs a})p(u, {\bs a})\ &{}=
\ (\tilde l - l)\,B(u, \bs h(\T))\,
\prod_{s=1}^n\prod_{j=1}^{m_s-1} (u-z_s+ j)\ .
\end{align*}
Hence, the polynomials $\tilde p(u, \tilde{\bs a}(\T))$\>,
\,$p(u, {\bs a}(\T))$ span a vector subspace $V$ lying in the intersection
$C_G$, see Theorem \ref{thm inverse on Schub}. Conversely, if $V$ is
a point of $ C_G$,
then $V$ has a basis $f,g$ like in Lemma \ref{lemma on diff eqn}.
Then by Lemma \ref{lemma on diff eqn} we have
\begin{align*}
\Wr (g(u),f(u))\ &{}=
\ (\tilde l - l)\,\prod_{s=1}^n\prod_{j=1}^{m_s} (u-z_s+ j)\ ,
\\
g(u)f(u-2) - g(u-2) f(u)\ &{}=
\ (\tilde l - l) \,\Bc(u)\,\prod_{s=1}^n\prod_{j=1}^{m_s-1} (u-z_s+ j)
\end{align*}
for a suitable polynomial $\Bc(u)$. Hence, the triple $g, f, \Bc$
determines a point $\T$, whose coordinates satisfy
equations \Ref{eqns for C_G}.
\end{proof}

\section{Algebras $A_P$ and $A_L$}
\subsection{Algebra $A_P$}
\label{sectAP}
Consider the space $\C^{\tilde l + l + n - 1}$
with coordinates $\tilde{\bs a}, \bs a, \bs h$. Let
\be
\Dh\ =\ d(u)\ - \ B(u, \bs h)\,\shift^{-1}\
+ \ a(u)\,\shift^{-2}\
\vvgood
\ee
be the difference operator defined in \Ref{D}. If $\bs h$ satisfies equations
$q_{1}(\bs h) = 0$ and $q_{2}(\bs h) = 0$, then the polynomial $\Dh(\qat)$ is
a polynomial in $u$ of degree $\tilde l+n-3$,
\vvn.2>
\be
\Dh(\qat)\ =\ \tilde {q}_3(\tilde{\bs a},\bs h)\,u^{\tilde l+n-3}
\ +\ \dots \ +\ \tilde q_{\tilde l+n}(\tilde{\bs a},\bs h)\ .
\vv.2>
\ee
The coefficients $\tilde{q}_i (\tilde{\bs a},\bs h)$ are functions linear
in $\tilde{\bs a}$ and linear in $\bs h$.

Recall that if $\qa = u^l + a_1u^{l-1}+ \dots + a_l$, and $\bs h$ satisfies
equations $q_{1}(\bs h) = 0$ and $q_{2}(\bs h) = 0$, then the polynomial
$\Dh(\qa)$ is a polynomial in $u$ of degree $l+n-3$,
\vvn.3>
\be
\Dh(\qa)\ =\ q_3(\bs a,\bs h)\,u^{l+n-3}\ +\ \dots \ +\
q_{l+n}(\bs a,\bs h)\ .
\vv.3>
\ee

Denote by $I_P$ the ideal in $\C[\tilde{\bs a}, \bs a, \bs h]$
generated by polynomials $q_1, q_2, q_3,\dots,q_{l+n}$,
$\tilde{q}_3,\dots,\alb\tilde q_{\tilde l+n}$.
The ideal $I_P$ defines a scheme $C_P\subset \C^{\tilde l + l + n - 1}$.
The algebra
\be
A_P\ =\ \C[\tilde{\bs a}, \bs a, \bs h]/ I_P\
\ee
is the algebra of functions on $C_P$.

The scheme $C_P$ is the scheme of points $\T\in \C^{\tilde l + l + n - 1}$
such that the difference equation $\Dhp w(u)=0$ has two polynomial solutions
$\tilde {p}(u,\tilde{\bs a}(\T))$ and $p(u,\bs a(\T))$.

\subsection{Isomorphism $\psi_{GP} : A_G \to A_P$ }
\subsubsection{}
\label{thm Isom psi-GP}
{\bf Theorem.}
{\it If $z_i-z_j \notin \Z$ for $i\ne j$, then the identity map
$\C^{\tilde l + l + n - 1} \to \C^{\tilde l + l + n - 1}$ induces
an algebra isomorphism $\psi_{GP} : A_G \to A_P$.
}

\begin{proof}
If $\T$ is a point of \,$C_P$, then polynomials
$\tilde p(u,\tilde{\bs a}(\T))$,
$p(u, {\bs a}(\T))$ are annihilated by the difference operator
$d(u) - B(u, \bs h(\T))\shift^{-1} + a(u) \shift^{-2}$.
If $z_i-z_j \notin \Z$ for $i\ne j$, then the span $V$ of polynomials
$\tilde p(u, \tilde{\bs a}(\T))$, $p(u,{\bs a}(\T))$ is a point of
\,$C_G$ by Theorem \ref{thm inverse on Schub}.
This reasoning defines an algebra homomorphism $\psi_{GP} : A_G \to A_P$.

Conversely, if $\T$ is a point of
$C_G$, then the triple $\tilde p(u, \tilde{\bs a}(\T))$,
$p(u, {\bs a}(\T))$, $B(u, \bs h(\T))$ satisfies equations
\vvn-.5>
\begin{align*}
\Wr\,(\tilde p(u, \tilde{\bs a}(\T)), p(u, {\bs a}(\T)))\ &{}=
\ (\tilde l - l)\,\prod_{s=1}^n\prod_{j=1}^{m_s} (u-z_s+ j)\ ,
\\
\tilde p(u, \tilde{\bs a}) p(u-2, {\bs a}) -
\tilde p(u-2, \tilde{\bs a})p(u, {\bs a})\ &{}=
\ (\tilde l-l)\,B(u,\bs h(\T))\,\prod_{s=1}^n\prod_{j=1}^{m_s-1}(u-z_s+ j)\ .
\end{align*}
Hence the polynomials $\tilde p(u, \tilde{\bs a}(\T))$,
$p(u, {\bs a}(\T))$ are annihilated by the difference operator
$d(u) - B(u, \bs h(\T))\shift^{-1} + a(u) \shift^{-2}$. Therefore, $\T$
is a point of $C_P$.
\end{proof}

\subsection{Algebra $A_L$}
\label{A_L}
Assume that $m_1,\dots,m_n,\,l$ satisfy conditions of
Section~\ref{sec new conditiong}.
Let $\bs \La =\alb(\La^{(1)}, \dots , \La^{(n)})$ be the collection of
$\glt$-highest weights with $\La^{(s)}=(m_s,0)$. Let
\vvn.2>
\be
L_{\bs \La}\ = \ L_{\La^{(1)}}\otimes \dots \otimes L_{\La^{(n)}}
\vv.2>
\ee
be the tensor product of irreducible $\glt$-modules with highest
weights $\La^{(1)}, \dots , \La^{(n)}$, respectively, and
$v_{\bs \La} = v_{(m_1,0)} \otimes \dots \otimes v_{(m_n,0)}$
the tensor product of the corresponding highest weight vectors.
Denote by
\be
L_{\bs \La}(\bs z)\ =
\ L_{\La^{(1)}}(z_1)\otimes \dots \otimes L_{\La^{(n)}}(z_n)
\vv.1>
\ee
the tensor product of evaluation modules.

\medskip
Let $\Snash \subset L_{\bs \La}(\bs z)$ be the subspace of $\glt$-singular
vectors of weight $(\sum_{i=1}^nm_i - l , l)$. The algebra $A_L$ is
the Bethe algebra associated with $\Snash$.

\medskip
Assume that $m_i\in\Z_{\geq 0}$ for
$i=1,\dots,n$, and $m_1\leq m_2\leq\dots\leq m_n$.
Assume that $z_i - z_j + m_j - s \ne 0$ and $z_i - z_j -1 - s \ne 0$
for all $i<j$ and $s=0,1,\dots,m_i-1$.
Then by Theorem \ref{sec product of irrep}, there is a natural isomorphism
$\W/K \to L_{\bs \La}(\bs z)$ such that $1 \mapsto v_{\bs \La}$.
Here $K\subset \W$ is the kernel of the Yangian Shapovalov form on $\W$.

The Yangian Shapovalov form on $\W$ induces the Yangian Shapovalov form
$S$ on $L_{\bs \La}(\bs z)$ such that $S(v_{\bs \La},v_{\bs \La})=1$
and $S(x\cdot v,\, w) = S(v, x^+\cdot w)$ for all $x\in \Y$ and
$v,w \in L_{\bs \La}(\bs z)$. The form $S$ is nondegenerate and symmetric.

We have the composition of linear maps
\vvn.1>
\be
\W\ \to\ \W/K\ \to\ L_{\bs \La}(\bs z)\ .
\vv.1>
\ee
Restricting this composition to $\Sing\,\W$ we get a linear epimorphism
\vvn.2>
\be
\sh \ :\ \snash \ \to\ \Snash\ .
\vv.2>
\ee
The Bethe algebra $A_W$ preserves the kernel of $\sh$ and induces
a commutative subalgebra in $\End\,(\Snash)$. The induced subalgebra coincides
with the Bethe algebra $A_L$. We denote by \,$\psi_{WL} : A_W \to A_L$
\,the corresponding epimorphism.

\medskip
The operators of the algebra $A_L$ are symmetric with respect
to the Yangian Shapovalov form on $L_{\bs \La}(\bs z)$.

\subsubsection{}
Denote by
\be
\D_L\ =\ d(u) \ - \
(2u^n + \psi_{WL}(H_1)u^{n-1} + \dots + \psi_{WL}(H_n))\,\shift^{-1}\ +
\ a(u)\,\shift^{-2}
\vv.2>
\ee
the universal difference operator associated with the subspace $\Snash$
and collection $\bs z$.

\subsubsection{}
\label{thm exist of polyn in Snash}
{\bf Theorem.}
{\it
Assume that the pair $\bs \La, l$\ satisfies conditions of
Section \ref{sec new conditiong}. Then for any $v_0\in \Snash$ there exist
$v_1,\dots,v_{\tilde l} \in \Snash$ such that the function
\vvn.2>
\be
w(u) \ =\ v_0\,u^{\tilde l} \, +\,v_1\,u^{\tilde l-1}\, +\, \dots \,+
\,v_{\tilde l}
\vv.2>
\ee
is a solution of the difference equation $\D_L w(u)\,=\,0$.
}

\medskip
This theorem is a particular case of Theorem 7.3 in \cite{MTV2}.

\section{Homomorphisms of algebras $A_D$, $A_P$ and $A_L$}
\label{New homomorphisms}

\subsection{Epimorphism $\psi_{DP} : A_D \to A_P$}
\label{sec Epi}

A point $\T$ of $C_P$ determines the difference equation $\Dhp\,w(u)\,=\,0$
and two solutions $\tilde {p}(u,\tilde{\bs a}(\T))$ and $p(u,\bs a(\T))$.
Then the pair, consisting of the difference operator $\Dhp$ and the solution
$p(u,\bs a(\T))$ of the smaller degree, determines a point of $C_D$,
see Section~\ref{subsub h's}. This correspondence defines a natural algebra
epimorphism $\psi_{DP} : A_D \to A_P$.

\subsection{Linear map $\xi : A_D \to \Snash$}
\label{Linear map xi}
Assume that $z_1,\dots,z_n, \ m_1,\dots, m_n$ satisfy the assumptions of
Theorem \ref{sec product of irrep}. Then we have the composition of linear maps
\vvn.2>
\be
A_D\ \stackrel{\phi}{\longrightarrow}\ A_D^*\
\stackrel{\tau}{\longrightarrow}
\ \snash\ \stackrel{\sh}{\longrightarrow}\ \Snash\ .
\vv.2>
\ee
Denote this composition by $\xi : A_D \to \Snash$.
By Theorem \ref{1st main thm}, \ $\xi$ is a linear epimorphism.

Let $\psi_{DL} : A_D \to A_L$ be the algebra epimorphism defined as the
composition $\psi_{WL}\>\psi_{DW}$.

\subsubsection{}
\label{lEmma on intertw}
{\bf Lemma.}
{\it
If $z_1,\dots,z_n, \ m_1,\dots, m_n$ satisfy the assumptions of
Theorem \ref{sec product of irrep}, then the linear map
$\xi$ intertwines the action of the multiplication operators $L_f,\,f\in A_D$,
on $A_D$ and the action of the Bethe algebra $A_L$ on $\Snash$, that is,
for any $f,g\in A_D$ we have\ $\xi(L_f(g)) \,=\, \psi_{DL}(f)(\xi(g))$.
}

\medskip
The lemma follows from Corollary \ref{intersting cor}.

\subsubsection{}
\label{tarasov's lem}
{\bf Lemma.}
{\it
If $z_1,\dots,z_n, \ m_1,\dots, m_n$ satisfy the assumptions of
Theorem \ref{sec product of irrep}, then the kernel of $\xi$ coincides with
the kernel of $\psi_{DL}$.
}

\begin{proof}
If $\psi_{DL}(f)=0$, then $\xi(f) = \xi(L_f(1)) = \psi_{DL}(f)(\xi(1)) = 0$.
On the other hand, if $\xi(f)=0$, then for any $g\in A_D$ we have
$\psi_{DL}(f)(\xi(g)) = \xi(L_f(g)) = \xi(fg) = \xi(L_g(f))
= \psi_{DL}(g)(\xi(f)) = 0$.
Since $\xi$ is an epimorphism, this means that $\psi_{DL}(f)=0$.
\end{proof}

\subsubsection{}
\label{2nd main lem}
{\bf Lemma.}
{\it
If $z_i-z_j \notin \Z$ for $i\ne j$, then the kernel of\, $\xi$ coincides
with the kernel of $\psi_{DP}$.
}

\begin{proof}
If $z_i-z_j \notin \Z$ for $i\ne j$, then the assumptions of
Theorem \ref{sec product of irrep} are satisfied and $\xi$ is defined.

By Schubert calculus,
\,$\dim\,\Snash\,=\,\dim\, A_G$. By Theorem \ref{thm Isom psi-GP}
$\dim\,A_G = \dim\,A_P$ if $z_i-z_j \notin \Z$ for $i\ne j$.
Hence it suffices to show that the kernel of $\xi$ contains the kernel
of $\psi_{DP}$. But this follows from Theorems \ref{thm exist of polyn}
and \ref{thm exist of polyn in Snash}.

Indeed the defining relations in $A_P = A_D/(\ker\,\psi_{DP})$ are the
conditions on the operator $\Dh$ to have two linear independent polynomials
in the kernel. Theorems \ref{thm exist of polyn} and
\ref{thm exist of polyn in Snash} guarantee these relations for elements
of the Bethe algebra $A_L$. Hence, the kernel of $\psi_{DL}$ contains
the kernel of $\psi_{DP}$. By Lemma \ref{tarasov's lem}, the kernel of $\xi$
coincides with the kernel of $\psi_{DL}$. Therefore, the kernel of $\xi$
contains the kernel of $\psi_{DP}$.
\end{proof}

\subsubsection{}
\label{Cor of 2nd main lem}
{\bf Corollary.}
{\it
Let $z_i-z_j \notin \Z$ for all $i\ne j$.
Then the algebras $A_P$, $A_L$ and $A_G$ are isomorphic.
}

\begin{proof}
Since the algebra epimorphisms $\psi_{DP}$ and $\psi_{DL}$ have
the same kernels, the algebras $A_P$ and $A_L$ are isomorphic.
Then $A_L$ and $A_G$ are isomorphic by Theorem~\ref{thm Isom psi-GP}.
\end{proof}

\medskip
\subsection{Second main theorem}
Let $z_i-z_j \notin \Z$ for all $i\ne j$.
Denote by $\psi_{PL} : A_P \to A_L$ the isomorphism induced by $\psi_{DL}$ and
$\psi_{DP}$. Lemmas~\ref{lEmma on intertw}\,--\,\ref{2nd main lem}
imply the following theorem.

\subsubsection{}
\label{second main thm}
{\bf Theorem.}
{\it
If $z_i-z_j \notin \Z$ for all
$i\ne j$, then the linear map $\xi$ induces a linear isomorphism
\vvn-.4>
\be
\zeta\ :\ A_P\ \to\ \Snash
\vv.3>
\ee
which intertwines the multiplication operators $L_f,\,f\in A_P$,\, on $A_P$
and the action of the Bethe algebra $A_L$ on $\Snash$, that is, for any
$f,g\in A_P$ we have $\zeta(L_f(g)) \,=\, \psi_{PL}(f)(\zeta(g))$.
}
\qed

\subsubsection{}
\label{Cor 3 of 2nd main thm}
{\bf Corollary. }
{\it
Let $z_i-z_j \notin \Z$ for all $i\ne j$. Assume that every operator
$f \in A_L$ is diagonalizable. Then the algebra $A_L$ has simple spectrum and
all points of the intersection of Schubert cycles
\vvn.1>
\be
C_G\ =\ C_{\F(\infty), \La^{(\infty)}} \,\tbigcap\,(\,\cap_{i=1}^n\,
C_{\F(z_i), \La^{(i)}}\,)\
\vv-.1>
\ee
are of multiplicity one.
}

\begin{proof}
The algebras $A_L$, $A_P$ and $A_G$ are isomorphic.
We have $A_P = \oplus_{\T} \,A_{\T,P}$ where the sum is over the points of the
scheme $C_P$ considered as a set and $A_{\T,P}$ is the local algebra associated
with a point $\T$. The algebra $A_{\T,P}$ has nonzero nilpotent elements
if $\dim\,A_{\T,P}>1$. If every element $f\in A_P$ is diagonalizable,
then the algebra $A_P$ is the direct sum of one-dimensional local algebras.
Hence $A_P$ has simple spectrum as well as the algebras $A_L$ and $A_G$.
\end{proof}

Corollary \ref{Cor 3 of 2nd main thm} has the following application.

\subsubsection{}
\label{cor on shapiro}
{\bf Corollary.}
{\it Assume that $\,z_1, \dots, z_n$ are real,
$z_i-z_j \notin \Z$ and $|z_i-z_j|\gg 1$ for all $i\ne j$.
Then all points of the intersection of Schubert cycles
\vvn.1>
\be
C_G\ =\ C_{\infty, \La^{(\infty)}} \,\tbigcap\,(\,\cap_{i=1}^n\,
C_{z_i, \La^{(i)}}\,)\
\vv-.1>
\ee
are of multiplicity one.}

\begin{proof}
If $z_1,\dots, z_n$ are real and $|z_i-z_j|\gg 1$ for all $i\ne j$,
then the Yangian Shapovalov form, restricted to the real part of
$\Sing\,L_{\bs \La}[\>l\>]$, is positive definite,
see Appendix~C in \cite{MTV1}. The Hamiltonians
$\psi_{WL}(H_1),$ \dots, $\psi_{WL}(H_1)$, restricted to the real part
of $\Sing\,L_{\bs \La}[\>l\>]$, are real symmetric operators operators
with respect to the Yangian Shapovalov form, see \cite{MTV1}.
Hence, all elements of the Bethe algebra
$A_L$ are diagonalizable operators. Therefore, the spectrum of $A_G$ is simple
and all points of $C_G$ are of multiplicity one.
\end{proof}

Corollary~\ref{cor on shapiro} is related to Theorem~1 from~\cite{EGSV}
and Theorem~2.1 from~\cite{MTV4} concerning the real Schubert calculus.

\subsubsection{}
{\bf Example.}
Let $n=3$, \,$\La^{(s)}=(1,0)$\>, \,${s=1,2,3}$, \,$\La^{(\infty)}=(2,1)$,
\,and
\be
R\,=\,4\>(z_1^2+z_2^2+z_3^2-z_1z_2-z_1z_3-z_2z_3)\>-\>3\;.
\ee
If $R\ne 0$, then every element of $A_L$ is diagonalizable and the algebra
$A_L$ is isomorphic to the direct sum $\C\oplus\C$. If $R=0$, then
the algebra $A_L$ contains a nonzero nilpotent matrix and is isomorphic
to \>$\C[b]/\langle b^2\rangle$.

\section{Operators with polynomial kernel and Bethe algebra $A_L$}
\label{Equations with polynomial solutions only}

\subsection{ Linear isomorphism $\theta : A_P^* \to \Snash$}
\label{Linear isomorphism theta}
Let $z_i-z_j \notin \Z$ for all $i\ne j$. Define the symmetric bilinear form
on $A_P$ by the formula \be (f,\,g)_P\ = \ S\bigl(\zeta(f), \,\zeta(g)\bigr)
\qquad {\rm for\ all} \quad f,g \in A_P \ , \ee where $S(\,,\,)$ denotes the
Yangian Shapovalov form on $\Snash$.

\subsubsection{}
{\bf Lemma.}
{\it
The form $(\,,\,)_P$ is nondegenerate.}

\medskip
The lemma follows from the fact that the Yangian Shapovalov form on $\Snash$
is nondegenerate and the fact that $\zeta$ is an isomorphism.

\subsubsection{}
{\bf Lemma.}
{\it
We have $(fg,h)_P = (g,fh)_P$ for all $f,g,h \in A_P$.}

\medskip
The lemma follows from the fact the elements of the Bethe algebra are symmetric operators
with respect to the Yangian Shapovalov form,
see Section~\ref{Algebra A_W and Shapovalov form}.

\medskip
The form $(\,,\,)_P$ defines a linear isomorphism $\pi:A_P \to A^*_P$,
$f\mapsto (f\,,\cdot)_P$.

\subsubsection{}
\label{COR}
{\bf Corollary.}
{\it
Let $z_i-z_j \notin \Z$ for all $i\ne j$. Then the map $\pi$
intertwines the multiplication operators $L_f,\,f\in A_P$, on $A_P$
and the dual operators $L^*_f,\,f\in A_P$,\, on $A^*_P$.
}

\subsection{Third main theorem}
Summarizing Theorem \ref{second main thm}
and Corollary \ref{COR}
we obtain the following theorem.

\subsubsection{}
\label{useful cor}
{\bf Theorem.}
{\it
Let $z_i-z_j \notin \Z$ for all $i\ne j$. Then the composition
$\theta=\zeta\pi^{-1}$ is a linear isomorphism from $A^*_P$ to $\Snash$
which intertwines the multiplication operators $L^*_f,\,f\in A_P$,
on $A^*_P$ and the action of the Bethe algebra $A_L$ on $\Snash$,
that is, for any $f\in A_P$ and $g\in A^*_P$ we have
$\theta(L^*_f(g))\,=\, \psi_{PL}(f)(\theta(g))$.}
\qed

\subsubsection{}
\label{Cor 2 of 2nd main thm NEW}
Let $z_i-z_j \notin \Z$ for all $i\ne j$.
Assume that $v\in \Snash$ is an eigenvector of the Bethe algebra $A_L$,
that is, $\psi_{WL}(H_s) v = \lambda_s v$ for suitable $\lambda_s\in\C$
and $s=1,\dots,n$. Then, by Corollary~7.4 in~\cite{MTV2},
the difference equation
\vvn.2>
\be
\bigl(d(u)\ -\ (2u^n + \la_1 u^{n-1} + \dots + \la_n)\,\shift^{-1} \ +
\ a(u)\,\shift^{-2}\bigr)\,w(u)\ =\ 0
\vv.2>
\ee
has two linearly independent polynomial solutions, one of degree $l$
and the other of degree $\tilde l$. The following corollary of
Theorem~\ref{useful cor} gives the converse statement.

\subsubsection{}
\label{Cor 2 of 2nd main thm}
{\bf Corollary of Theorem \ref{useful cor}.}
{\it
Let $z_i-z_j \notin \Z$ for all $i\ne j$.
Assume that $(\la_1,\dots,\la_n)\in \C^n$
is a point such that
\be
\la_1\ =\ \sum^n _{i=1}\,(m_i-2z_i) \, ,
\qquad
\la_2\ =\ l\,\Bigl(l - 1 -\sum_{i=1}^n\,m_i\Bigr)\ +
\sum_{1\leq i< j\leq n}\!\bigl(z_iz_j + (z_i-m_i)(z_j-m_j)\bigr)\,,
\ee
and the difference equation
\vvn.2>
\beq
\label{diffeq}
\bigl(d(u)\ -\ (2u^n + \la_1 u^{n-1} + \dots + \la_n)\,\shift^{-1} \ +
\ a(u)\,\shift^{-2}\bigr)\,w(u)\ =\ 0
\vv.2>
\eeq
has two linearly independent polynomial solutions. Then there exists
a unique up to normalization eigenvector $v \in \Snash$ of the action of
the Bethe algebra $A_L$ such that for every $s=1,\dots,n$ we have}
\vvn.2>
\beq
\label{Hv=lav}
\psi_{WL}(H_s) \,v \,=\, \la_s\, v\ .
\vv.3>
\eeq

\begin{proof}[Proof of Corollary \ref{Cor 2 of 2nd main thm}]
\ Indeed, such a point $(\la_1,\dots,\la_n)$ defines a linear function
${\eta : A_P \to \C}$, \,$h_s \mapsto \la_s$, for $s=1,\dots,n$. Moreover,
\,$\eta(fg) = \eta(f)\eta(g)$ for all $f,g \in A_P$. Hence $\eta\in A_P^*$ is
an eigenvector of operators $L^*_f$ acting on $A_P^*$.
By Theorem~\ref{useful cor}, the vector $v=\theta(\eta)\in\Snash$
is an eigenvector of the action of the Bethe algebra $A_L$ with eigenvalues
prescribed in Corollary~\ref{Cor 2 of 2nd main thm}.

Let $v'\in\Snash$ satisfy~\Ref{Hv=lav}, then $\eta'=\theta^{-1}(v)\in A_P^*$
satisfies \,$\eta'(fg)=\eta(f)\eta'(g)$ for all $f,g \in A_P$. Hence, for $g=1$
we have \,$\eta'(f)=\eta(f)\eta'(1)$. Therefore, $\eta'$ is proportional to
$\eta$, and $v'$ is proportional to $v$.
\end{proof}

\subsubsection{}
\label{how to find}
Assume that $(\la_1,\dots,\la_n)\in\C^n$ is a point satisfying the assumptions
of Corollary~\ref{Cor 2 of 2nd main thm}. We describe how to find the
eigenvector $v\in\Snash$, indicated in Corollary~\ref{Cor 2 of 2nd main thm}.

Let $f(u)$ be the monic polynomial of degree $l$ which is a solution of
the difference equation~\Ref{diffeq}. Consider the polynomial
\be
\om(\bs y)\ = \ y_0^l \prod_{j=1}^{n-1} f(y_j-1)
\ee
as an element of $\Wy$, see Section \ref{Universal weight function}.
By Theorem \ref{thm on Bethe ansatz} this vector lies in
$\Sing\,\Wy[\>l\>]$ and $\om(\bs y) $ is an eigenvector of the Bethe algebra
$A_W$ with eigenvalues prescribed in Corollary~\ref{Cor 2 of 2nd main thm}.
Consider a maximal subspace $V \subset \snash$ with three properties:
\\ \phantom{aa}
{}\ i) \,$V$ contains $\om(\bs y)$,
\\ \phantom{aa}
{}\ ii) \,$V$ does not contain other eigenvectors of the Bethe algebra $A_W$,
\\ \phantom{aa}
{}\
iii) \,$V$ is invariant with respect to the Bethe algebra $A_W$.
\\
Such a maximal subspace does exist and is unique. Let $\sh(V)\subset\Snash$
be the image of $V$ under the epimorphism $\sh$. Then by
Corollary~\ref{Cor 2 of 2nd main thm}, the subspace $\sh(V)$ contains
a unique one-dimensional subspace of eigenvectors of the Bethe algebra $A_L$.
Any such an eigenvector may serve as an eigenvector of the Bethe algebra $A_L$
indicated in Corollary~\ref{Cor 2 of 2nd main thm}.

\section{Homogeneous \XXX/ Heisenberg model}
\label{Heisenberg chain}

\subsection{Statement of results}
\label{Statement of results}
In Sections \ref{Three more algebras}--\ref{Equations with polynomial solutions
only}, in most of the assertions we assumed that $z_1,\dots,z_n \in \C$ are
such that $z_i-z_j\notin\Z$ for $i\ne j$, and $m_1,\dots,m_n$ are natural
numbers. In this section we assume that \beq
\label{01}
z_1=\dots =z_n =0 \qquad {\rm and} \qquad m_1=\dots =m_n =1\ .
\eeq
This special case is called {\it the homogeneous \XXX/ Heisenberg model}.

In other words, in this section we consider the $\Y$-module
\be
L_{\bs 1}(\bs 0)\ =\ L_{(1,0)}(0) \otimes \dots \otimes L_{(1,0)}(0)\ ,
\ee
which is the tensor product of $n$ copies of the two-dimensional evaluation
module, and the subspace of $\glt$-singular vectors of weight $(n-l,l)$,
\vvn.2>
\be
\Sing\,L_{\bs 1}[\>l\>]\ =\
\{\, p \in L_{\bs 1}(\bs 0)\ | \ e_{12} p = 0, \ e_{22} p = l p\, \}\ .
\vv.2>
\ee
The subspace $\SSnash$ is not empty if and only if \>$2l \le n$, that is,
if and only if
the pair $\bigl((1,0),\dots,\alb (1,0)\bigr)\>,\,l$ is separating. In that case
\vvn-.4>
\be
\dim\,\Sing\,L_{\bs 1}[\>l\>]\ = \ \binom nl-\binom n{l-1}\ .
\ee

\medskip
The algebra $A_L$ is the Bethe algebra associated with the subspace
$\Sing\,L_{\bs 1}[\>l\>]$. It is generated by the coefficients of the series
$\bigl(T_{11}(u)+T_{22}(u)\bigr)\big|_{\SSnash}$\,.

\medskip
The main result of this section is the following theorem.

\subsubsection{}
\label{thm on Heisenberg}
{\bf Theorem.}
{\it
For the homogeneous \XXX/ Heisenberg model,
the Bethe algebra $A_L$ has simple spectrum.}

\medskip
The theorem will be proved in Section \ref{proofs-H}.

\subsubsection{}
\label{thm on Wronskian}
Let \,$\tilde l = n+1-l$. We have \,$\tilde l+l-1=n$ and \,$\tilde l>l$.
Denote by $f,g$ \,two polynomials in $\C[u]$ of the form:
\begin{align}
\label{two polyns}
f(u)\ &{}=\ u^l + f_1u^{l-1}+ \dots + f_l\ ,
\\
g(u)\ &{}=\ u^{\tilde l} + {g}_1u^{\tilde l-1}+ \dots +
g_{\tilde l - l -1} u^{l+1} +
g_{\tilde l - l +1}u^{l-1}+ \dots +g_{\tilde l}\ .
\notag
\end{align}

As a byproduct of the proof of Theorem \ref{thm on Heisenberg}
we prove the following theorem.

\medskip
\noindent
{\bf Theorem.}
{\it
There exist exactly \,$\binom nl-\binom n{l-1}$ distinct pairs of polynomials
$f,g$ \>of the form \Ref{two polyns}, such that\/}
\be
f(u)g\>(u-1)-f(u-1)\>g(u)\ =\ (l-\tilde l)\>(u+1)^n\ .
\ee

\medskip
Theorem \ref{thm on Wronskian} will be proved in Section \ref{proofs-W}.

\subsection{Algebra $A_L$ for the homogeneous \XXX/ model}
Consider the Yangian module $\W$ corresponding to the polynomials
\be
a(u)\ =\ (u+1)^n\ ,\qquad d(u)\ =\ u^n\ .
\ee
The numbers \Ref{01} satisfy the assumptions of
Theorem \ref{sec product of irrep}. Therefore
the $\Y$-module $L_{\bs 1}(\bs 0)$ is
irreducible, and there is a natural epimorphism $\W \to L_{\bs 1}(\bs 0)$
of $\Y$-modules. Restricting this epimorphism to $\Sing\,\W[\>l\>]$,
\vvn.3>
we obtain a linear epimorphism
\be
\sh \ :\ \snash \ \to\ \Snash\ .
\vv.3>
\ee
The Bethe algebra $A_W$ preserves the kernel of $\sh$ and induces
a commutative subalgebra in $\End\,(\Snash)$. The induced subalgebra coincides
with the Bethe algebra $A_L$, see~Section~\ref{A_L}.

Denote by \,$\psi_{WL} : A_W \to A_L$ \,the corresponding epimorphism. We have
\vvn.2>
\be
\bigl(T_{11}(u) + T_{2}(u)\bigr)\big|_{\,\SSnash}\ = \
2 + \psi_{WL}(H_1)\>u^{-1} + \dots + \>\psi_{WL}(H_n)\>u^{-n} \ ,
\ee
where
\vvn-.4>
\be
\psi_{WL}(H_1)\ =\ n\ , \qquad\psi_{WL}(H_2)\ =\ l (l-1-n) + \frac{n(n-1)}2\ ,
\ee
see Section \ref{H-Hamiltonians}. Thus the Bethe algebra $A_L$ is generated
by elements $\psi_{WL}(H_3), \dots,$ $\psi_{WL}(H_n)$.

\subsection{Algebra $A_P$ for the homogeneous \XXX/ model}
Consider the space $\C^{2n}$ with coordinates $\tilde{\bs a}, \bs a, \bs h$,
as in Section \ref{sec Presentation of algebra A_G}, and polynomials
$\qat$, $\qa$, $B(u, \bs h)$.

Given the polynomials $a(u)=(u+1)^n$ and $d(u)=u^n$, we define the ideal
$I_P$, the algebra $A_P$, and the scheme $C_P$ as in Section~\ref{sectAP}.
The scheme $C_P$ is the scheme of points $\T\in \C^{2n}$ such that
the difference equation
\be
\label{D heis}
(u^n\ - \ B(u, \bs h)\,\shift^{-1}\ + \ (u+1)^n\,\shift^{-2})\,w(u)=0
\vv.3>
\ee
has two polynomial solutions $\tilde p(u,\tilde{\bs a}(\T))$ and
$p(u,\bs a(\T))$.

\subsection{Algebra $A_G$ for the homogeneous \XXX/ model}
\label{A_G heis}
Consider the space $\C^{2n}$ with coordinates $\tilde{\bs a}, \bs a, \bs h$,
and polynomials $\qat$, $\qa$, $B(u, \bs h)$. Let us write
\begin{align*}
& \Wr\,(\qat, \qa)\ = \ (\tilde l - l) u^{n} +
w_1(\tilde {\bs a}, \bs a)u^{n-1} + \dots + w_n(\tilde {\bs a}, \bs a)\ ,
\\[4pt]
& \tilde p(u, \tilde{\bs a}) p(u-2, {\bs a}) -
\tilde p(u-2, \tilde{\bs a})p(u, {\bs a}) \ ={}
\\[3pt]
& \hphantom{\Wr\,(\qat, \qa)\ ={}\!}
2 (\tilde l - l) u^{n} +
\hat w_1(\tilde {\bs a}, \bs a)u^{n-1}+\dots+\hat w_n(\tilde{\bs a},\bs a)
\end{align*}
for suitable polynomials $w_1,\dots,w_n$, $\hat w_1,\dots, \hat w_n$
in variables $\tilde{\bs a}, \bs a$.

Denote by $I_G$ the ideal in $\C[\tilde{\bs a}, \bs a, \bs h]$ generated by
$2n$ polynomials
\beq
\label{eqns for C_G - heis}
w_i(\tilde{\bs a}, \bs a) - (\tilde l - l) \binom{n}{i}\ ,\qquad
\hat w_i(\tilde{\bs a}, \bs a) - (\tilde l - l) h_i\ ,\qquad
i = 1, \dots , n\ .
\eeq
The ideal $I_G$ defines a scheme $C_G \subset \C^{2n}$. Then
\be
A_G\ =\ \C[\tilde{\bs a}, \bs a, \bs h]/ I_G
\ee
is the algebra of functions on $C_G$.

The scheme $C_G$ is the scheme of points $\T \in \C^{2n}$ such that
\vvn.3>
\begin{gather}
\label{GC1}
\Wr\, (\tilde p (u, \tilde{\bs a}(\T), p (u, {\bs a}(\T))\;=
\;(\tilde l - l)\,(u+1)^n\ ,
\\[3pt]
\tilde p(u, \tilde{\bs a}) p(u-2, {\bs a}) -
\tilde p(u-2, \tilde{\bs a})p(u, {\bs a}) \;=
\;(\tilde l - l)\, B(u,\bs h(\T))\ .
\notag
\end{gather}

\subsubsection{}
\label{thm Isom psi-GP heis}
{\bf Theorem.}
{\it The identity map $\C^{2n} \to \C^{2n}$ induces an algebra isomorphism
$\psi_{GP} : A_G \to A_P$.
}

\begin{proof}
The proof is similar to the proof of Theorem \ref{thm Isom psi-GP}.
\end{proof}

\subsubsection{}
{\bf Lemma.}
\label{sec dim A_G H}
{\it The dimension of $A_G$ considered as a vector space is equal to}
\be
\dim\,\Sing\,L_{\bs 1}[\>l\>]\ = \ \binom nl-\binom n{l-1}\ .
\vvn.2>
\ee

\begin{proof}
Consider the ideal $I_G(\bs z)$ defined by \Ref{eqns for C_G}
for $m_1=\dots=m_n=1$ and arbitrary $z_1,\dots,z_n$. Consider the algebra
$A_G(\bs z)\,=\,\C[\tilde{\bs a}, \bs a, \bs h]/ I_G(\bs z)$.
By Lemma \ref{lemma on new presentation}, if $z_1,\dots,z_n$ are
distinct and close to zero, then $A_G(\bs z)$ is the algebra of
functions on the intersection of Schubert cells $C_G(\bs z)$, see~\Ref{CG},
and by \Ref{dimAG} we have
\be
\dim\,A_G(\bs z)\ =\ \binom nl-\binom n{l-1}\ .
\ee
To complete the proof of Lemma \ref{sec dim A_G H}, it suffices to verify
two facts:
\begin{enumerate}
\item[(i)]
There are no algebraic curves over $\C$ lying in the scheme $C_G(\bs 0)$,
defined by the ideal \Ref{eqns for C_G - heis}.

\item[(ii)]
Let a sequence $\bs z^{(i)}$, $i=1,2,\dots{}$, tend to $\bs0$.
Let $\T^{(i)} \in C_G(\bs z^{(i)})$\>, \,$i = 1, 2, \dots{}\,$,
be a sequence of points. Then all coordinates
$\bigl(\bs{\tilde a}(\T^{(i)}), \bs a(\T^{(i)}), \bs h(\T^{(i)}\bigr)$
remain bounded as $i$ tends to infinity.
\end{enumerate}

By Theorem~\ref{thm Isom psi-GP}, the schemes $C_G(\bs z)$ and $C_P(\bs z)$
are isomorphic if $z_1,\dots,z_n$ are distinct and close to zero.
By Theorem~\ref{thm Isom psi-GP heis}, the schemes $C_G(\bs 0)$ and
$C_P(\bs 0)$ are isomorphic as well. Claims (i) and (ii) hold for the scheme
$C_P(\bs z)$ by Theorem~\ref{thm const wrt z} because $C_P(\bs z)$ is
a subscheme of the scheme $C_D(\bs z)$.
\end{proof}

\subsection{Three more homomorphisms for the homogeneous \XXX/ model}
\label{sec EPI-h}
In Sections \ref{sec Epi} and \ref{Linear map xi}, we define
an algebra epimorphism $\psi_{DP} : A_D \to A_P$, \,a linear epimorphism
$\xi : A_D \to \SSnash$ as the composition of linear maps
\be
A_D\ \stackrel{\phi}{\longrightarrow}\
A_D^*\ \stackrel{\tau}{\longrightarrow}
\ \snash\ \stackrel{\sh}{\longrightarrow}\ \SSnash\ .
\ee
and an algebra epimorphism $\psi_{DL} : A_D \to A_L$ as the composition
$\psi_{WL}\psi_{DW}$.

For the homogeneous \XXX/ model, we have Lemmas~\ref{lEmma on intertw}
and~\ref{tarasov's lem} and the following
analogue of Lemma \ref{2nd main lem}.

\subsubsection{}
\label{2nd main lem H}
{\bf Lemma.}
{\it For the homogeneous \XXX/ model, the kernel of\, $\xi$ coincides
with the kernel of $\psi_{DP}$.
}

\begin{proof}
The proof is similar to the proof of Lemma~\ref{2nd main lem} with
Theorem~\ref{thm Isom psi-GP heis} replacing Theorem~\ref{thm Isom psi-GP}.
\end{proof}

\subsubsection{}
\label{Cor of 2nd main lem h}
{\bf Corollary.}
{\it
For the homogeneous \XXX/ model, the algebras $A_P$, $A_L$ and $A_G$ are
isomorphic.
}

\medskip
Denote by $\psi_{PL} : A_P \to A_L$ the isomorphism induced by $\psi_{DL}$ and
$\psi_{DP}$. We have the following analogue of Theorem~\ref{second main thm}.

\subsubsection{}
\label{second main thm h}
{\bf Theorem.}
{\it
For the homogeneous \XXX/ model, the linear map $\xi$ induces a linear
isomorphism
\vvn-.3>
\be
\zeta\ :\ A_P\ \to\ \SSnash
\vv.2>
\ee
which intertwines the multiplication operators $L_f,\,f\in A_P$,\, on $A_P$
and the action of the Bethe algebra $A_L$ on $\SSnash$, that is, for any
$f,g\in A_P$ we have $\zeta(L_f(g)) \,=\, \psi_{PL}(f)(\zeta(g))$.
}
\qed

\subsection{The Bethe algebra $A_L$ of the homogeneous \XXX/ model
is diagonalizable}
\subsubsection{}
\vsk-.2>
\label{thm on diag}
{\bf Theorem.}
{\it
For the homogeneous \XXX/ model, all elements of $A_L$ are diagonalizable
operators.}

\begin{proof}
Let $v_+$ be a highest $\glt$-weight vector of $L_{(1,0)}$ and
$v_- = e_{21}v_+$. Then $v_+,v_-$ form a basis of $L_{(1,0)}$.
Consider the Hermitian form on $L_{\bs 1}(\bs 0)$ for which the vectors
\bea
v_{i_1} \otimes \dots \otimes v_{i_n}
\qquad
{\rm with}
\qquad
i_j \in \{+,-\}
\eea
generate an orthonormal basis of $L_{\bs 1}(\bs 0)$. For any
$X\in\End\bigl(L_{\bs 1}(\bs 0)\bigr)$, denote by $X^\dag$ the Hermitian
conjugate operator with respect to this Hermitian form. It is clear that
\vvn.2>
\be
\bigl((1^{\otimes(j-1)}\otimes e_{ab}\otimes 1^{\otimes(n-j)})
|_{L_{\bs 1}(\bs 0)}\bigr)^{\dag}\ =
\ (1^{\otimes(j-1)}\otimes e_{ba}\otimes 1^{\otimes(n-j)})
|_{L_{\bs 1}(\bs 0)}\ .
\vv.2>
\ee

Using the fact that $(e_{11}+e_{22})|_{L_{(1,0)}}=1$ and the definition
of the coproduct~\Ref{Delta}, it is straightforward to verify by induction
on $n$ that
\vvn.3>
\be
\bigl(T_{ab}(u)|_{L_{\bs 1}(\bs 0)}\bigr)^{\dag}\ =
\ (-1)^{a+b+n}\,T_{3-a,3-b}(-\>\bar u-1)|_{L_{\bs 1}(\bs 0)}\ ,
\vv.3>
\ee
where $\bar u$ is the complex conjugate of $u$. Therefore,
\vvn.3>
\be
\bigl((T_{11}(u)+T_{22}(u))|_{L_{\bs 1}(\bs 0)}\bigr)^{\dag}\ =
\ -\>(T_{11}(-\>\bar u-1)+T_{22}(-\>\bar u-1))|_{L_{\bs 1}(\bs 0)}\ .
\vvgood
\vv.2>
\ee
This means that for any $X\in A_L$, the Hermitian conjugate operator $X^\dag$
lies in $A_L$. Hence, any element of $A_L$ commutes with its Hermitian
conjugate and, therefore, is diagonalizable.
\end{proof}

\subsection{Proof of Theorem \ref{thm on Heisenberg}}
\label{proofs-H}
The proof is similar to the proof of Corollary~\ref{Cor 3 of 2nd main thm},
because every element of $A_L$ is diagonalizable by Theorem~\ref{thm on diag}.
\qed

\subsection{Proof of Theorem \ref{thm on Wronskian}}
\label{proofs-W}
The algebras $A_G$ and $A_L$ are isomorphic. So, by Theorem~\ref{thm on diag}
every element $f\in A_G$ is diagonalizable. Therefore, the algebra $A_G$ is
the direct sum of one-dimensional local algebras. Hence $C_G$ considered as
a set consists of \,$\dim\,A_G\binom {n}{l}-\binom{n}{l-1}$ distinct points,
see Lemma~\ref{sec dim A_G H}. Theorem \ref{thm on Wronskian} is proved.

\subsubsection{}
\label{Cor 2 of 2nd main thm NEW h}
Assume that $v\in \SSnash$ is an eigenvector of the Bethe algebra $A_L$,
that is, $\psi_{WL}(H_s) v = \lambda_s v$ for suitable $\lambda_s\in\C$
and $s=1,\dots,n$. Then by Corollary~7.4 in~\cite{MTV2},
the difference equation
\vvn.2>
\be
\bigl(u^n \ -\ (2u^n + \la_1 u^{n-1} + \dots + \la_n)\,\shift^{-1} \ +\
(u+1)^n\,\shift^{-2}\bigr)\,w(u)\ =\ 0
\vv.2>
\ee
has two linearly independent polynomial solutions, one of degree $l$
and the other of degree $n-l+1$. The following corollary of
Theorem \ref{thm on Heisenberg} gives the converse statement.

\subsubsection{}
\label{Cor 2 of 2nd main thm h}
{\bf Corollary of Theorem \ref{thm on Heisenberg}.}
{\it
Assume that \,$(\la_1,\dots,\la_n)\in \C^n$ is a point such that
\vvn-.6>
\be
\la_1\ =\ n\ ,
\qquad \la_2\ =\ l(l-1-n) - \frac{n(n-1)}2\ ,
\ee
and the difference equation
\be
\bigl(u^n \ -\ (2u^n + \la_1 u^{n-1} + \dots + \la_n)\,\shift^{-1} \ +\
(u+1)^n\,\shift^{-2}\bigr)\,w(u)\ =\ 0
\ee
has two linearly independent polynomial solutions. Then there exists
a unique up to normalization eigenvector $v \in \SSnash$ of the action of
the Bethe algebra $A_L$ of the homogeneous \>\XXX/ \<model such that for every
$s=1,\dots,n$ we have}
\be
\psi_{WL}(H_s) \,v \,=\, \la_s\, v\ .
\ee

The proof of Corollary~\ref{Cor 2 of 2nd main thm h}] is similar
to the proof of Corollary~\ref{Cor 2 of 2nd main thm}.

\subsubsection{}
\label{how to find h}
Assume that $(\la_1,\dots,\la_n)\in\C^n$ is a point satisfying the assumptions
of Corollary~\ref{Cor 2 of 2nd main thm h}. In order to find the eigenvector
$v\in\SSnash$, indicated in Corollary \ref{Cor 2 of 2nd main thm h}, one needs
to apply the procedure described in Section \ref{how to find}.

\end{document}